\documentclass{nm}
\renewcommand\thisnumber{x}
\renewcommand\thisyear {201x}

\renewcommand\thisvolume{xx}
\usepackage{charter}
\usepackage[charter]{mathdesign}
\usepackage{amsmath}
\usepackage{graphicx}%
\usepackage{latexsym ,rawfonts}
 \usepackage{multirow}
\usepackage{tikz}
\usetikzlibrary{shapes}
\usepackage{color}
\usepackage{caption}
\usepackage{subfigure}
\usepackage{algorithmic}
\usepackage{savesym}
\savesymbol{iint}
\usepackage{txfonts}
\restoresymbol{TXF}{iint}

\usepackage{caption}
\usepackage{subfigure}
\usepackage{algorithmic}

\def \bbf{{f}}
\def \bg{{g}}

\def \bu{{ u}}

\def \bw{{ w}}
\def \bx{{ x}}

\def \bH{{ H}}

\def \bL{{ L}}

\def \bX{{ X}}

\def \bsigma{{\boldsymbol \sigma}}

\def \dt{{\Delta t}}

\def \det{\mathrm{det}}

\def \bubble{\mathrm{bubble}}
\def \bouss{\mathrm{bouss}}
\def \phinew{\phi^{\mathrm{new}}_h}
\def \psinew{\psi^{\mathrm{new}}_h}
\def \phiold{\phi^{\mathrm{old}}_h}
\def \psiold{\psi^{\mathrm{old}}_h}

\def \eps{\epsilon}
\def \veps{\varepsilon}

\def \sE{\mathscr{E}}
\def \sF{\mathscr{F}}

\def \wtilde{\widetilde}
\def \pa{\partial}

\begin{document}

\markboth{J. Brannick, C. Liu, T. Qian, and H. Sun}{Diffuse Interface Methods for Multiple Phase Materials}
\title{Diffuse Interface Methods for Multiple Phase Materials: An Energetic Variational Approach}
\author[J. Brannick, C. Liu, T. Qian, and H. Sun]{J. Brannick \affil{1}\comma\corrauth, C. Liu\affil{1}, T. Qian\affil{2}, and H. Sun\affil{1}}
\address{\affilnum{1}\ Department of Mathematics, The
       Pennsylvania State University, University Park, PA 16802, USA.
 \\
 \affilnum{2}\ Hong Kong University of Science and Technology, Clear Water Bay, Kowloon
Hong Kong.
       }

\emails{
 {\tt brannick@psu.edu} (J. Brannick),
  {\tt liu@math.psu.edu} (C. Liu),
 {\tt tmaqian@ust.hk} (T. Qian),
 {\tt sun@math.psu.edu} (H. Sun)
 }

\begin{abstract}
In this paper, we introduce a diffuse interface model for describing
the dynamics of mixtures involving multiple (two or more) phases.
The coupled hydrodynamical system 
 is derived through an energetic variational approach.
The total energy of the system includes the kinetic energy and
the mixing (interfacial) energies. The least action principle (or
the principle of virtual work) is applied to derive the conservative
part of the dynamics, with a focus on the reversible part of
the stress tensor arising from the mixing energies. The dissipative
part of the dynamics is then introduced through a dissipation function
in the energy law, in line with the Onsager principle of
least energy dissipation. The final system, formed by a set of coupled
time-dependent partial differential equations, reflects a balance
among various conservative and dissipative forces and governs
the evolution of velocity and phase fields. To demonstrate
the applicability of the proposed model, a few two-dimensional
simulations have been carried out, including (1) the force balance
at the three-phase contact  line in equilibrium, (2) a rising bubble penetrating
a fluid-fluid interface, and (3) a solid particle falling in
a binary fluid. The effects of slip at solid surface have been examined
in connection with contact line motion and a pinch-off phenomenon.
\end{abstract}




\maketitle



\section{Introduction}
Phase field methods (PFM), also known as diffuse interface methods,  have been widely used in 
modeling two-phase problems and free interface motion of mixtures.
The methods  are based on a labeling function $\phi(x)$, which usually takes values as
$\pm 1$,  to distinguish between the two different materials (phases).   
Du et. al. applied phase field methods in their studies of the configurations and the deformations of elastic bio-membranes
\cite{DuLiRhWa05a}. Liu and Shen investigated the use of two-phase models for studying bubble relaxation,
rise, and coalescence \cite{LiSh02}. Qian et al. studied 
the moving contact line problem using phase field methods 
in \cite{QWS03}. Yue et. al. \cite{YuFeLiSh05b}
studied a general approach for modeling two-phase complex fluids, with 
numerical examples simulating emulsion of nematic drops in a Newtonian
matrix. Recently Shen and Yang applied the phase field method to
two-phase incompressible flows with different densities and viscosities \cite{SY10}.

The basic idea of the two-phase PFM is to use a coarse graining (mean field) model to describe the 
microscopic dynamics of the mixtures.  In the hydrodynamical (macroscopic) time scale,
such dynamics involve the deformations of each phase,
the  interaction between the two, and their interactions with 
the surrounding environment.
The  underlying dynamical system is derived from applying variational principles  to  a certain free energy,
e.g. the classical Ginzburg-Landau type energy \cite{CaHi58}
\[
\sF_{CH} = 
\int \gamma \Big\{ \,\frac{\veps}{2} \, |\nabla \phi|^2 +
\frac{1}{4\veps} \, \big( \phi^2 -1 \big)^2 \, \Big\} dx , 
\]
where $\phi(x)$ is the phase field function and $\veps$ is the
width of the diffuse interface. The two  parts in the above integrand represent
the ``philic" and ``phobic" interactions between the two materials. 
The parameter $\gamma$ can be associated to the surface tension in
the conventional sharp interface formulations.
The applicability of this model has been demonstrated for many different
applications  (see \cite{DuLiRyWa05} \cite{DuLiWa04} \cite{QiWaSh06}
and references theirin).
Although analytically it is still an open question whether the
the sharp interface model can be recovered by the phase field model via a
rigorous proof, the latter has been applied theoretically and
numerically for a long time. Moreover, from a practical and more physical point of view,
the sharp interface models can be viewed as the simplification  or idealization 
of phase field models. 

\begin{figure}[ht!]
 \centerline{ \includegraphics[width=.4\textwidth]{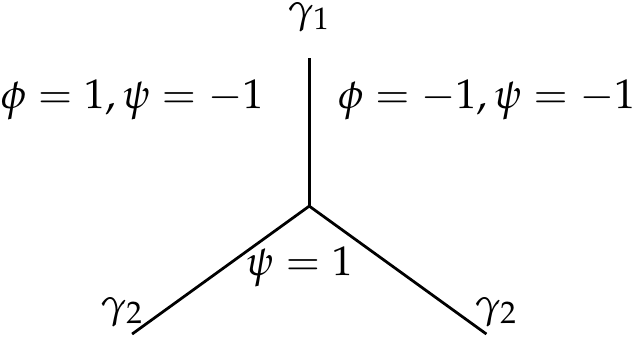}}
  \caption{
    A Schematic illustration 
    for three phases
    distinguished/labelled by two phase fields $\phi$ and $\psi$.
    \label{fig:threephase}}
\end{figure}

In this paper, we show that for problems in which 
more than two phases are involved, we can 
introduce additional labeling functions to distinguish among them,
as illustrated in \textsc{Figure} \ref{fig:threephase}.  The derivation follows from applying 
the energetic variational framework as
in \cite{YuFeLiSh05b, LiSh02}. 
Here, in the region at the bottom of the figure, a single phase is characterized by $\{\psi = 1\}$ and 
$\phi$ is not defined.  In the top region of the figure, there are two phases distinguished
by different values of $\phi$, while sharing the same $\psi$
value. In a similar way, four different phases  can be characterized by two phase field 
functions.  We note that such an approach has been considered in other contexts \cite{BGS96}.

The remaining sections of the paper are organized as follows. We derive our multi-phase model using a variational approach in Section~\ref{sec:EVA}.
In Section~\ref{sec:APPS}, we discuss the numerical methods
in our simulations and present results for a variety of 2d multi-phase applications.  
 Section~\ref{sec:Conclude} presents  concluding remarks and future  work.

\section{Derivation of multi-phase model}\label{sec:EVA}
In analogy to approaches used for modeling two-phase problems~\cite{LiSh02}, 
for multi-phase problems we define a mixing energy
{\small
\begin{eqnarray}\label{eqn: elastic energy W}
E = \int_\Omega W(\phi, \nabla \phi, \psi, \nabla \psi) dx &=&
 \int_\Omega \biggr\{ \gamma_1 \, \Big( \frac{\psi -1}{2} \Big)^2
\Big( \frac{\veps_1}{2}|\nabla \phi|^2 + \frac{1}{4\veps_1}(\phi^2
-1)^2 \Big)    \\\nonumber
 +
 \quad  \gamma_2\Big(\frac{\veps_2}{2} |\nabla \psi|^2 &+& 
\frac{1}{4\veps_2} (\psi^2 - 1)^2 \Big)
  \biggr\} dx , 
\end{eqnarray}}
\noindent
where $\veps_1$ and $\veps_2$ are the widths of the interfaces along 
the differing phases and $\gamma_1$ and $\gamma_2$ are the surface tensions 
at the interfaces.
The $(\frac{\psi-1}{2})^2$ coefficient in the energy density for
$\phi$ is used to ensure that interactions between two different phases
(labeled as $\phi =  1,  \psi = -1$ and $\phi = -1, \psi = -1$) do not
directly influence the bulk of the third phase (labeled as $\psi = 1$).
We assume both $\gamma_1$ and $\gamma_2$
are constants, although in general
the parameters $\veps_i$, $\gamma_i$, $i=1,2$ can be
taken to be phase-dependent. 

Here we want to stress the relation of our approach to other related treatments by
other groups, such as those by Kim and Lowengrub \cite{Ki07,LoKi05}. 
It is clear, one needs at least 2 phase field functions to label the 3 distinct materials.
The key difference, and hence one of the main difficulties, is in the choices of the free energy functionals. 
While the free energy functionals in \cite{Ki07,LoKi05} are nondegenerate,
the free energy in \eqref{eqn: elastic energy W} involves degeneracy for the region 
$\{x: \psi = 1\}$. This generic degeneracy stands for the physics that in this solid
region, there is no effects from the fluids interactions (due to the dynamics of $\phi$).

Next, by adding the fluid equations to the system,  we obtain the total
energy of the hydrodynamic system as a weighted sum
of the kinetic energy and the mixing energy
\[
\sE = \int_\Omega \, \Big( \, \frac{1}{2} \rho \, |{ \bu}|^2 +
\lambda\, W(\phi, \nabla \phi, \psi, \nabla \psi) \, \Big) dx.
\]
Here, the constant $\lambda$ measures the competition between the two
types of energy.
To derive the stress tensor from the Ginzburg-Landau energy we apply the 
principle of virtual work (PoVW) \cite{DoEd86}, which states that the virtual work of the 
elastic energy $$E =  \int_\Omega W(\phi, \nabla \phi, \psi, \nabla \psi)$$ due to a virtual displacement $\delta x$ is
given by
\begin{equation} \label{eq:PoVW}
  \delta E = \int_\Omega \sigma : \nabla \delta \bx dx 
  = - \int_\Omega    
  (\nabla \cdot \sigma ) \cdot \delta \bx  dx.
\end{equation}
Before calculating the left-hand side in (\ref{eq:PoVW}), we 
take the macroscopic kinematic assumption of both  $\phi$ and $\psi$ being
 convective only, i.e,  there is no relaxation which contributesto the 
microscopic internal dissipation:
\begin{eqnarray} 
  &&\phi_t + \bu \cdot \nabla \phi = 0, \qquad
  \psi_t + \bu \cdot \nabla \psi = 0.\label{eq:pure convection}\\
  &\Rightarrow&
  \delta \phi + \delta \bx \cdot \nabla \phi = 0, \qquad
  \delta \psi + \delta \bx \cdot \nabla \psi = 0.
  \nonumber
\end{eqnarray}
Taking the gradient of the two equations in (\ref{eq:pure convection}) gives 
\begin{eqnarray*}
  \delta \nabla \phi + \left( \nabla \delta \bx \right)^T \nabla \phi 
  + (D^2 \phi) \delta \bx &=&0,\\
  \delta \nabla \psi + \left( \nabla \delta \bx \right)^T \nabla \psi 
  + (D^2 \psi) \delta \bx &=&0,
\end{eqnarray*}
where $D^2 \phi$ is the Hessian  of $\phi$. 

The virtual work can then be calculated as follows:
{\small
\begin{eqnarray*}
  \delta E &=&  
   \int_\Omega\left( \frac{\pa W}{\pa \phi} \cdot \delta \phi +  
  \frac{\pa W}{\pa \nabla \phi} : \delta \nabla \phi +
  \frac{\pa W}{\pa \psi} \cdot \delta \psi  +  
  \frac{\pa W}{\pa \nabla \psi} : \delta \nabla \psi \right)dx \\
  &=& \int_\Omega \gamma_1\left(\frac{\psi-1}{2}\right)^2
  \frac{(\phi^2 - 1) \phi}{\veps_1^2}   
  \left( - \delta \bx \cdot \nabla \phi \right) dx 
  +
  \int_\Omega  \gamma_1\left(\frac{\psi-1}{2}\right)^2
  \left[
    -(\nabla \delta \bx)^T \nabla \phi - \left(D^2 \phi \right)
    \, \delta \bx
    \right] dx \\
  && + \int_\Omega \gamma_2\left\{ \frac{(\psi^2 - 1)\psi}{\veps_2^2}
  + \frac{\psi-1}{2} \left[\frac{1}{2}|\nabla \phi|^2 
    + \frac{1}{4 \veps_1}\left(\phi^2-1 \right)\right] \right\} 
  \left( - \delta \bx \cdot \nabla \psi \right)dx \\
  && + \int_\Omega \gamma_2 \nabla \psi \, 
  \left[ - \left( \nabla \delta \bx \right)^T \,\nabla \psi
     - \left( D^2 \psi \right) \, \delta \bx     \right] dx \\
  &=& - \int_\Omega \gamma_1 \left(\frac{\psi-1}{2} \right)^2 
  \nabla \frac{\left(\phi^2-1\right)^2}{4 \veps_1^2} \, \delta \bx 
  - \int_\Omega \gamma_1\left[ \left(\frac{\psi-1}{2}\right)^2
    \nabla \phi\otimes \nabla \phi\right] \,\nabla \delta \bx dx \\
 && - \int_\Omega\gamma_1 \left(\frac{\psi-1}{2} \right)^2 \nabla
  \frac{|\nabla \phi|^2}{2} \, \delta \bx dx
   -\int_\Omega \nabla \frac{(\psi^2 -1)^2}{4 \veps_2^2} \,\delta
  \bx dx  \\
  && -\int_\Omega \gamma_2 \nabla \left[\left(\frac{\psi-1}{2}
    \right)^2\right] \,\left[ 
    \frac{1}{2} |\nabla \phi|^2 + \frac{1}{4 \veps_1} 
    (\phi^2 - 1)^2 
    \right] \, \delta \bx dx \\
  && - \int_\Omega \gamma_2\left(\nabla \psi \otimes
  \nabla \psi\right) \,\nabla \delta \bx dx 
   - \int_\Omega \gamma_2 \nabla \frac{|\nabla \psi|^2}{2} \,
   \delta \bx dx.
\end{eqnarray*}
}
Here, $\nabla$ denotes the gradient with respect to $\bx$.

Another integration by parts then gives
\begin{eqnarray*}
  \delta E &=& \int_\Omega\left[ -\nabla 
    W(\phi, \nabla \phi, \psi, \nabla \psi)
    \cdot \delta \bx - \left(\frac{\pa W}{\pa \nabla \phi}
    \otimes \nabla \phi
    + \frac{\pa W}{\pa \nabla \psi}
    \otimes \nabla \psi
    \right) : \nabla \delta \bx
        \right] dx \\
  &=& \int_\Omega -\nabla \cdot
  \left( 
  W \, I - \frac{\pa W}{\pa \nabla \phi}
    \otimes \nabla \phi
    - \frac{\pa W}{\pa \nabla \psi}
    \otimes \nabla \psi
  \right) \delta \bx  dx, 
\end{eqnarray*}
where $I$ is the identity matrix. 
The above gives the induced elastic force due to the interfacial mixing
energy.
Assuming the mixture of
incompressible fluids such that $\nabla \cdot \delta \bx = 0$,
the elastic stress tensor is uniquely determined up to 
an isotropic stress tensor, $f(\bx) I$. 
Therefore, the elastic stress tensor 
due to the mixing energy becomes
\[
\bsigma^e =   W \, I - \frac{\pa W}{\pa \nabla \phi}
    \otimes \nabla \phi
    - \frac{\pa W}{\pa \nabla \psi}
    \otimes \nabla \psi
    + f(\bx) I .
\]
We note that both isotropic tensors above can be absorbed into the pressure 
gradient term, which leads to the simplified equation 
\[
\wtilde\bsigma^e =    - \frac{\pa W}{\pa \nabla \phi}
    \otimes \nabla \phi
    - \frac{\pa W}{\pa \nabla \psi}
    \otimes \nabla \psi 
    = - \gamma_1 \veps_1 \left(\frac{\psi-1}{2} \right)^2 \nabla \phi \otimes
    \nabla \phi - \gamma_2 \veps_2 \nabla \psi \otimes \nabla \psi.
\]
We omit this dependence on tilde in the remainder of the paper, i.e.,
we use $\bsigma^e$ for the simplified stress tensor. We mention 
that a more precise definition of incompressibility is given by
$J = \det \frac{\pa \bx}{\pa\bX} = 1$.  Further, we note that a variation with respect to the domain can be
used to account for this constraint, as 
discussed in \cite{SL08}.


Next, we add the dissipative terms to the system, namely the viscous stress
tensor
$\bsigma^{v} = \frac{\nabla \bu + (\nabla \bu)^T}{2}$ with 
the viscosity coefficient $\mu$.
We also introduce the dissipation terms 
$\frac{\delta E}{\delta \phi}$ and $\frac{\delta E}{\delta \psi}$
into the convection relaxation equations for $\phi$ and $\psi$,
with $M_1$ and $M_2$ denoting the rate coefficients
for the relaxation at the interfaces.
The resulting system of PDEs is now as follows
\begin{eqnarray}\label{eq:MPM}
  \rho\, ( \bu_t + \bu\cdot \nabla \bu ) + \nabla p
  &=& \lambda \, \nabla \cdot (\bsigma^e + \bsigma^v ) \label{eqn: NSE} \: ,\\
  \phi_t + \bu \cdot \nabla \phi & = & - M_1 \frac{\delta
    E}{\delta \phi} \label{eqn: phi} \: , \\
  \psi_t + \bu \cdot \nabla \psi & = & - M_2 \frac{\delta E}{\delta
    \psi} \: ,  \label{eqn: psi}
\end{eqnarray}
where
\begin{eqnarray*}
  \frac{\delta E}{\delta \phi} &=& -\gamma_1 \left\{ 
\veps_1 \nabla \cdot \left[\left( \frac{\psi - 1}{2} \right)^2 \nabla \phi
  \right] - \frac{1}{\veps_1} \left( \frac{\psi -1}{2} \right)^2
\left( \phi^2 - 1\right) \phi
  \right\} \: , \\
\end{eqnarray*}
and 
\begin{eqnarray*}
  \frac{\delta E}{\delta \psi} &=& -\gamma_2 \left\{
  \veps_2 \Delta \psi - \frac{1}{\veps_2} (\psi^2 - 1) \psi
  \right\} + \gamma_1 \veps_1\frac{\psi - 1}{2} \left\{
  \frac{1}{2}|\nabla \phi|^2 + \frac{1}{4 \veps_1} (\phi^2 - 1)^2
  \right\}.
\end{eqnarray*}
\noindent
The above convection-relaxation equations (\ref{eqn: phi}) and (\ref{eqn:
psi}) can be interpreted as a fastest decent method for the energy. Recall that the
coefficients $M_1$ and $M_2$ determine the rates of the relaxation. 
The total system is thus dissipative with the governing energy law derived
by multiplying (\ref{eqn: NSE}) by $\bu$, (\ref{eqn: phi}) by
$\frac{\delta E}{\delta \phi}$ and (\ref{eqn: psi}) by $\frac{\delta
  E}{\delta \psi}$, integrating and adding the results together, and
then integrating by parts once again:
\begin{equation} \label{eqn: energy law}
\frac{d\sE}{dt}  = -\int_\Omega \biggr( \mu |\nabla \bu|^2 +
{\lambda}{M_1} \Big|\frac{\delta E}{\delta \phi} \Big|^2 +
{\lambda}{ M_2} \Big|\frac{\delta E}{\delta \psi}\Big|^2 \biggr)
\,dx.
\end{equation}

An {\it a priori} estimate of the solution then follows from 
the energy law:
\begin{eqnarray*}
\bu \in L^\infty\left( 0, T, \bL^2(\Omega) \right) \cap L^2 \left(
0, T, \bH^1(\Omega)\right)\\
\psi \in L^\infty \left(0,T, \bH^1(\Omega) \right) \cap L^2 \left(
0,T, \bH^2(\Omega) \right)
\end{eqnarray*}
However, the regularity of $\phi$ can not be determined from
the energy law because of the (possibly) degenerate pre-factor 
$(\psi -1)^2/4$. 

We understand that the relaxational equations (\ref{eqn: phi}) and
(\ref{eqn: psi}) do not
lead to the conservation of the order parameters $\phi$ and $\psi$.
This issue can be solved by using the Cahn-Hilliard dynamics or
a Lagrangian multiplier. As the main purpose of this paper is to
introduce a mixing free energy for the three-phase mixture problem,
we would like to leave the issue of order parameter conservation
to our future works. In our numerical simulations,
the rate coefficients $M_1$ and $M_2$ have been carefully chosen
to control the violation of conservation.

\begin{remark} In simulations one can ``turn off" the fluid by
setting the velocity to zero in which case the motion is purely driven
by mean curvature \cite{DuLiRhWa05a}.  In similar ways, we can add 
other mechanisms to adapt to various models.
\end{remark}

\begin{remark}  We mention that we can also set
$M_1$ and $M_2$ to zero to reflect the pure
transport kinematics.  Analytical results on the zero Weissenberg
number case can be found in \cite{LiLiZh05} \cite{LeLiZh08}.
Although this case is of great interest in certain applications, simulation of the pure transport equation is beyond the scope of the current paper.
\end{remark}

\section{Multi-phase simulations}\label{sec:APPS}
In this section, we use  numerical simulations to illustrate the applicability of our PFM for 
various multi-phase models. 

\subsection{Force Balance}
To begin, we consider three different phases of materials as shown in 
Figure~\ref{fig:threephase}.
In the absence of external forces and fluid effects, the dynamics of the system 
are driven purely by forces due to surface tensions induced from the 
mixing energy $E[\phi, \psi]$.
An imbalance of forces at the three-phase contact line will thus drive the
morphology until the forces are balanced.
For example, if the surface tensions on the interfaces are
equal, i.e. $\gamma_1 = \gamma_2$,, then at equilibrium all three angles
formed by the phases at the junction point are equal ($120^\circ$).

To model this system numerically, we omit the fluid equations in \eqref{eq:MPM} 
and use piecewise linear finite elements for the equations for $\phi, \psi$ spatially.  
Temporally, we use a fully implicit newton iteration.   
The discretized system
is
\begin{eqnarray*}
  \frac{\phinew - \phi^n_h}{\dt} &=& M_1 \gamma_1 
  \nabla \cdot \left[ \left( \frac{\psi_h^n-1}{2}\right)^2 
    \nabla \phinew \right] \\ && 
  - M_1 \gamma_1 \frac{1}{\veps_1^2} \left( \frac{\psi_h^n-1}{2}\right)^2
  \left[ ((\phiold)^2-1)\phiold
  + (3 (\phiold)^2 -1) (\phinew - \phiold)
  \right]\: , 
  \\
  \frac{\psinew - \psi^n_h}{\dt} &=& M_2 \gamma_2 
  \left[ \Delta \psinew - \frac{1}{\veps_1^2} ((\psiold)^2 - 1)\psiold
    -\frac{1}{\veps_1^2} (3 (\psiold)^2 - 1) (\psinew- \psiold)
    \right]\\
  &&- M_2 \gamma_1 \left( \frac{\psinew - 1}{2} \right)
  \left[ \frac{1}{2} |\nabla \phi_h^n|^2 +  
    \frac{1}{4 \veps_1^2 }((\phi_h^n)^2-1)^2
    \right] \: .
\end{eqnarray*} 
Here, $\dt$ is the
time step; $\phi_h^n$ and $\psi_h^n$ are the  finite element
solutions at previous  time steps $t_n$; and 
$\phinew$, $\phiold$, $\psinew$ and $\psiold$ are the Newton iterates 
at time step $t_{n+1}$.
The discrete equations lead to symmetric yet possibly indefinite linear
systems to solve at each time step for both $\phi$ and $\psi$.  Our choice
of solver is an ILU(0)-preconditioned GMRES method~\cite{EES83}.
To track the moving interfaces throughout a simulation we use adaptive 
mesh coarsening and refinement~\cite{HuReRu94} 
based on the Kelly error estimator \cite{KGZB83}, defined as
follows
\[
e^2 = \sum_{i=1}^{N^h}\frac{h}{24} \int_{\pa \Omega_i^h} J^2 ds, 
\]
where $J$ is the jump across the element boundary in the finite
element approximation to the gradient. Here, $\{\Omega_i^h\}_{i=1}^{N_h}$
is the partition of the computational domain.
In the simulation we compute the Kelly error estimator for the
linear combinations
of $\phi$ and $\psi$ and refine those elements with the largest estimated
errors that together make up 80 percent of the error and coarsen those cells 
of the error that account for a combined 10 percent of the smallest error 
(see Figure \ref{fig:bubble} for an illustration). 

\begin{figure}[ht!]
\centerline{
  \subfigure[\label{fig:phaseAngle}]{\includegraphics[width=.33\textwidth]{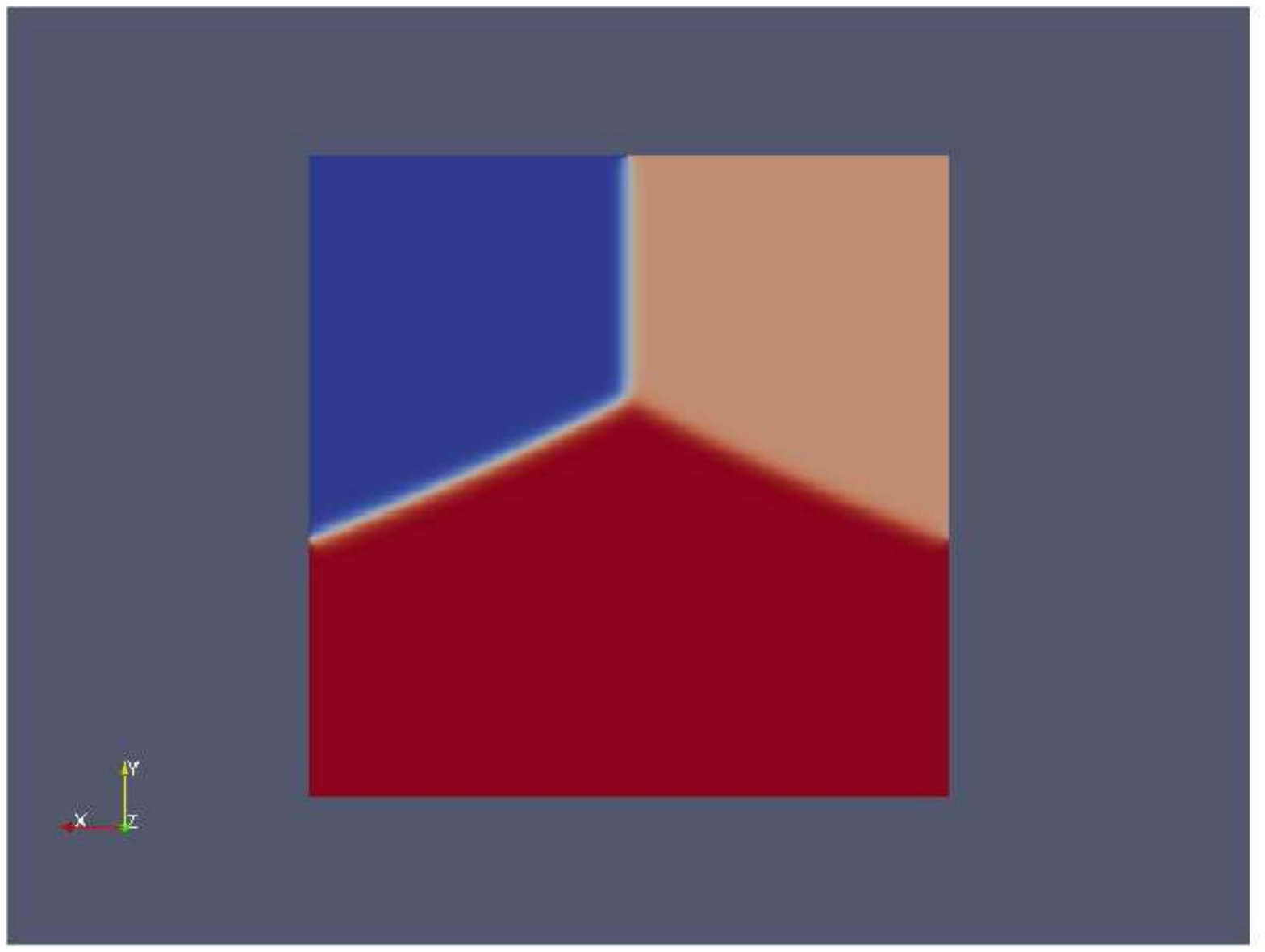}}
  \hfill
  \subfigure[\label{fig:phaseAngle1}]{\includegraphics[width=.33\textwidth]{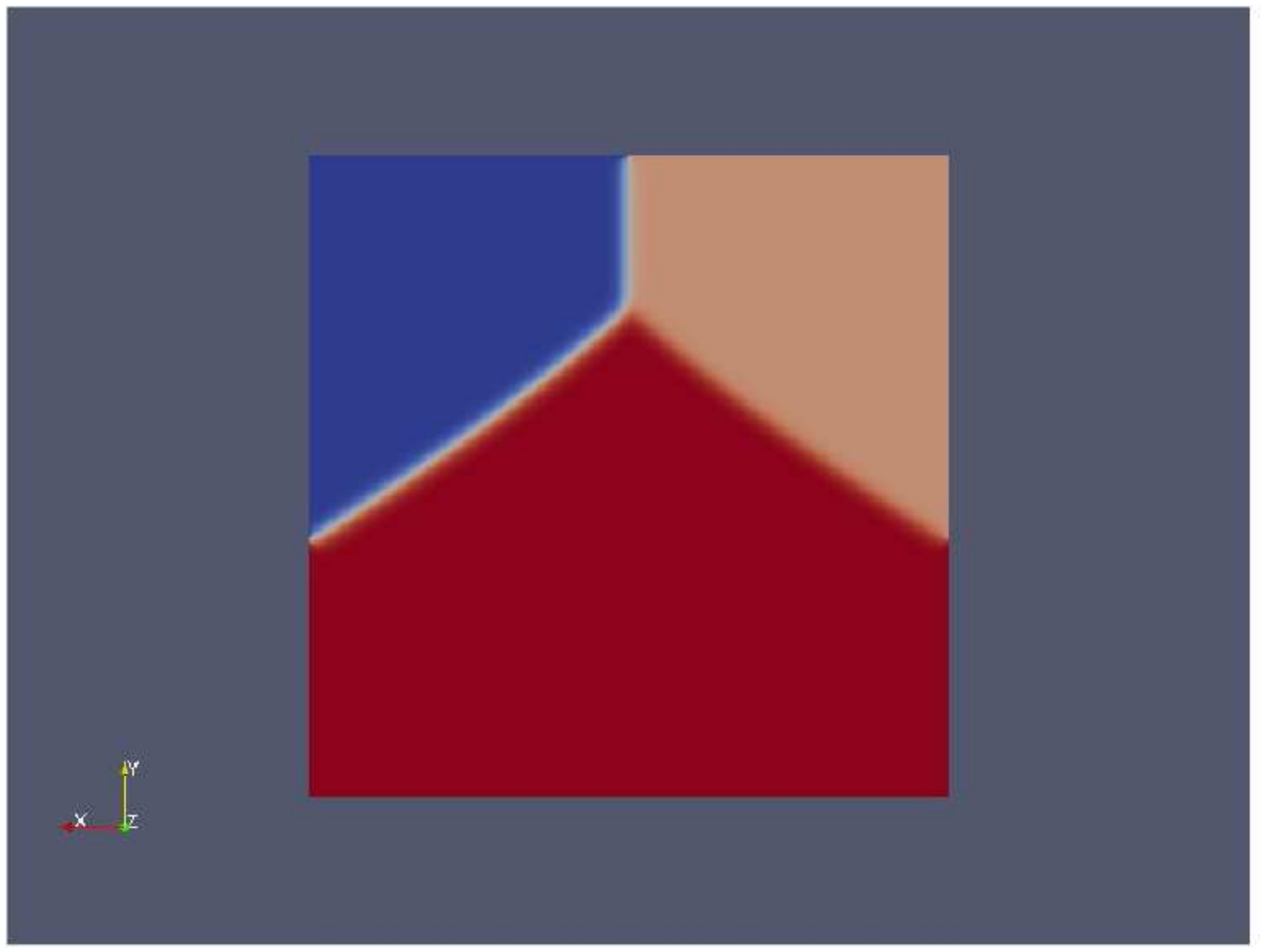}}
  \hfill
  \subfigure[\label{fig:phaseAngle2}]{\includegraphics[width=.33\textwidth]{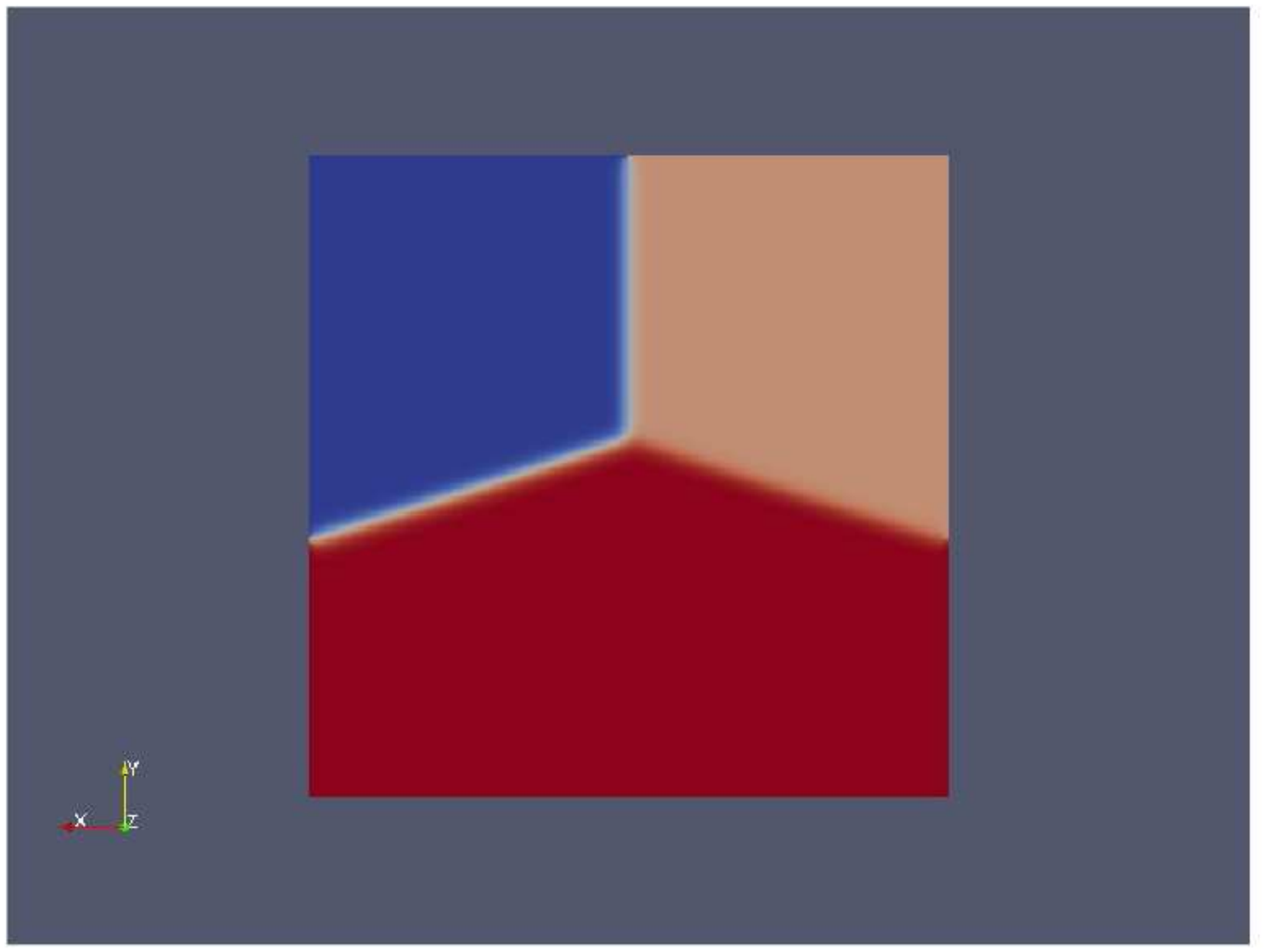}}
}
  \caption{ 
    Force balance at the three-phase contact line in equilibrium
    for different choices of surface tensions $\gamma_1$ and $\gamma_2$:
    \ref{fig:phaseAngle}: $\gamma_1 = \gamma_2 = 1.0$;  
    \ref{fig:phaseAngle1}: $\gamma_1 = 1.5$, and  $\gamma_2 = 1.0$; 
    \ref{fig:phaseAngle2}: $\gamma_1= 1.0$, and  $\gamma_2=1.5$. 
    \label{fig:ForceBalance}
    In each case, the three dihedral angles and
    the three surface tensions ($\gamma_1$, $\gamma_2$, and $\gamma_2$)
    form a Neumann's triangle.
  }
\end{figure}

The computational domain in this simulation is $[0,1]\times[0,
1]$, with the parameters set as $\mu = 1.0$, 
 $\dt = 0.1$, $\veps_1 =\veps_2 = 0.01$,  and
$M_1 = M_2 = 0.001$. The surface tensions $(\gamma_1, \gamma_2)$ 
are chosen as 
$(1.0, 1.0)$, $(1.5, 1.0)$ and $(1.0, 1.5)$ for three simulation
cases.
The initial conditions for the simulation are set as
\[
\phi(x,y) = \left\{
\begin{array}{cc}
  -1 & x < 0.5 \\
  +1 & x > 0.5
\end{array}
\right. \quad \mbox{and} \quad
\psi(x,y) = \left\{
\begin{array}{cc}
  -1 & y < 0.4 \\
  +1 & y > 0.4
\end{array}
\right. .
\]
The results of three numerical tests for various choices
of the surface tensions are provided in Figure~\ref{fig:ForceBalance}.
We note that, as expected, the angles between phases
at the three-phase contact line increase as the surface tension decreases
relative to its value on the other interfaces.

\subsection{A Rising Bubble}
The second application we consider is a model of a rising fluid bubble penetrating
an interface formed by two different fluids. Here we label the
fluid bubble by $\{ \psi = 1 \}$ and the region outside of it by $\{\psi = -1\}$.
The two other fluids are then distinguished by an additional phase-field 
labeling function $\phi$. We increase the densities of the outside
fluids to induce
a buoyancy which moves the bubble upward. Further, we mention that
to avoid the complications that arise from variable
densities, we adopt the Boussinesq approximation \cite{LiSh02}.
Additionally, we assume a low Reynolds number for the flow and thus 
replace the Navier-Stokes equations with the Stokes
equations.  
With these simplifications, the multi-phase model reduces to
\begin{eqnarray}
&&-\mu \Delta \bu + \nabla P= -\lambda \,\,\nabla \cdot \biggr\{
\widetilde\gamma_1\,\Big( \,
 \frac{\psi-1}{2} \Big)^2 \, \nabla \phi \otimes \nabla \phi
+ \widetilde\gamma_2 \nabla \psi \otimes \nabla \psi \biggr\}  + \bbf_{\bouss} \: , 
\label{nse bubble}\\
&&\nabla \cdot \bu = 0 \: , \\
&&\bbf_{\bouss}  =   -(1+\psi) \, (\rho_\bubble - \rho_0) \, \bg -
(1-\psi) \left[(1+\phi)\,(\rho_1 - \rho_0) + (1-\phi) \, (\rho_2 -
\rho_0) \right]\bg \: . \label{eqn: boussinesq}
\end{eqnarray}
Here, $\rho_0$ is the background density and its difference from
the actual densities ($\rho_\bubble$, $\rho_1$, or $\rho_2$) 
gives rise to the buoyancy force, with $\rho_\bubble$ denoting 
the density of the fluid in the bubble and $\rho_1$ and $\rho_2$
denoting the densities of the other two fluids.
The resulting system involves equations \eqref{eqn: phi}, \eqref{eqn: psi}, and
\eqref{nse bubble} - \eqref{eqn: boussinesq}.
\begin{figure}[ht!]
\begin{center}
\subfigure[\label{fig:bubble0}]{\includegraphics[width=.17\textwidth]
    {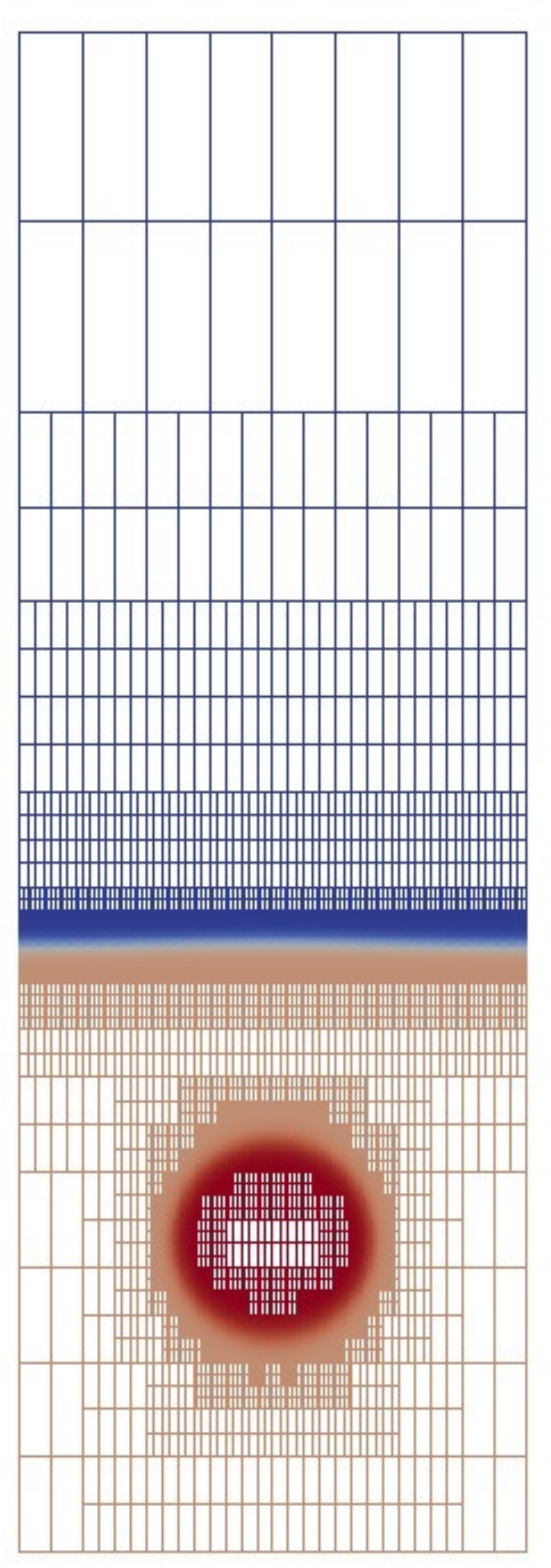}}
  \subfigure[\label{fig:bubble1}]{\includegraphics[width=.17\textwidth]
    {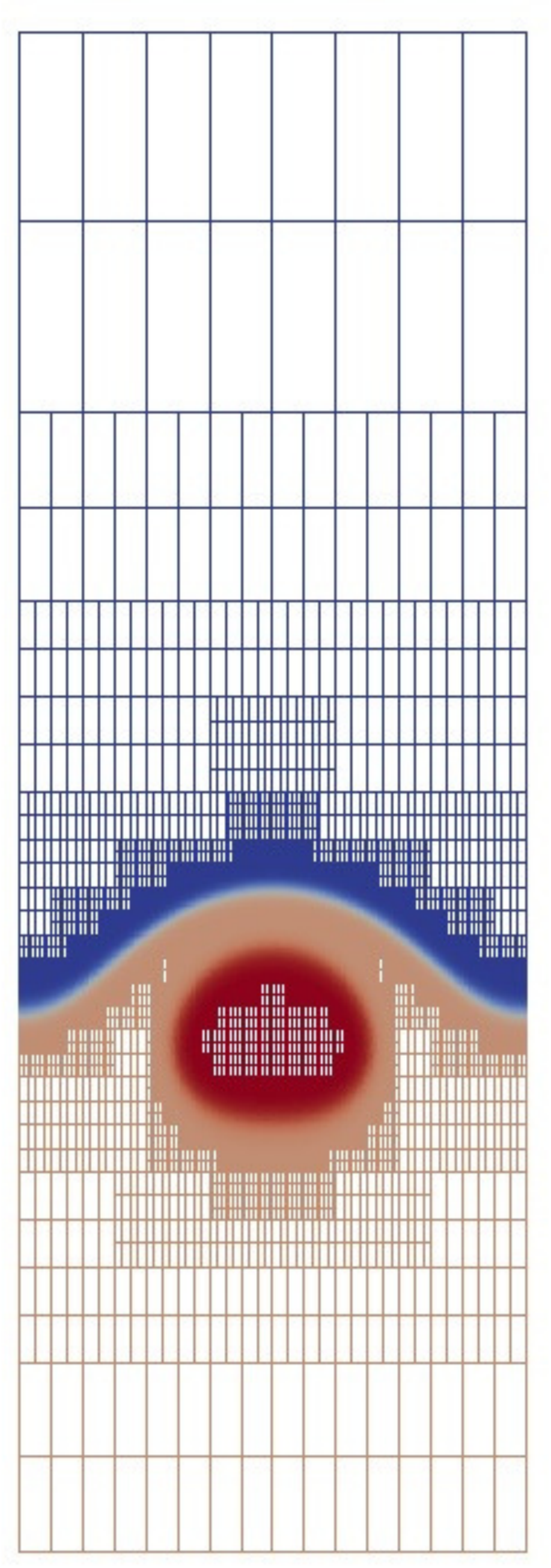}}
  \subfigure[\label{fig:bubble2}]{\includegraphics[width=.17\textwidth]
    {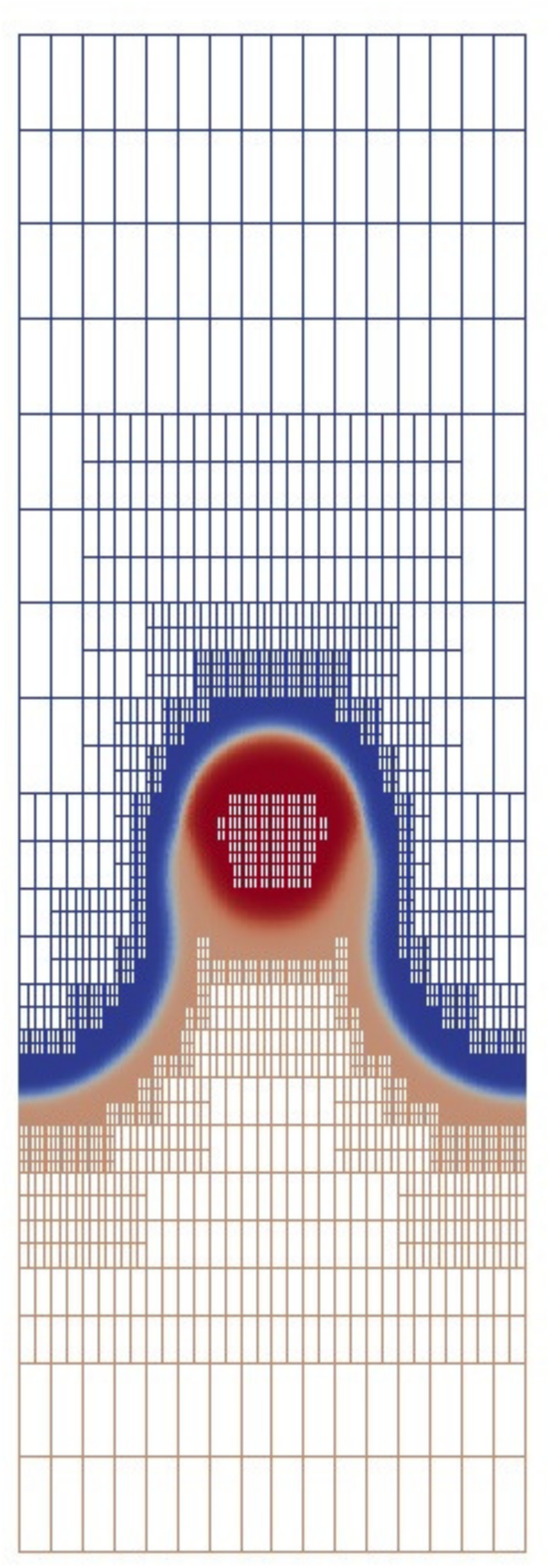}} 
  \subfigure[\label{fig:bubble3}]{\includegraphics[width=.17\textwidth]
    {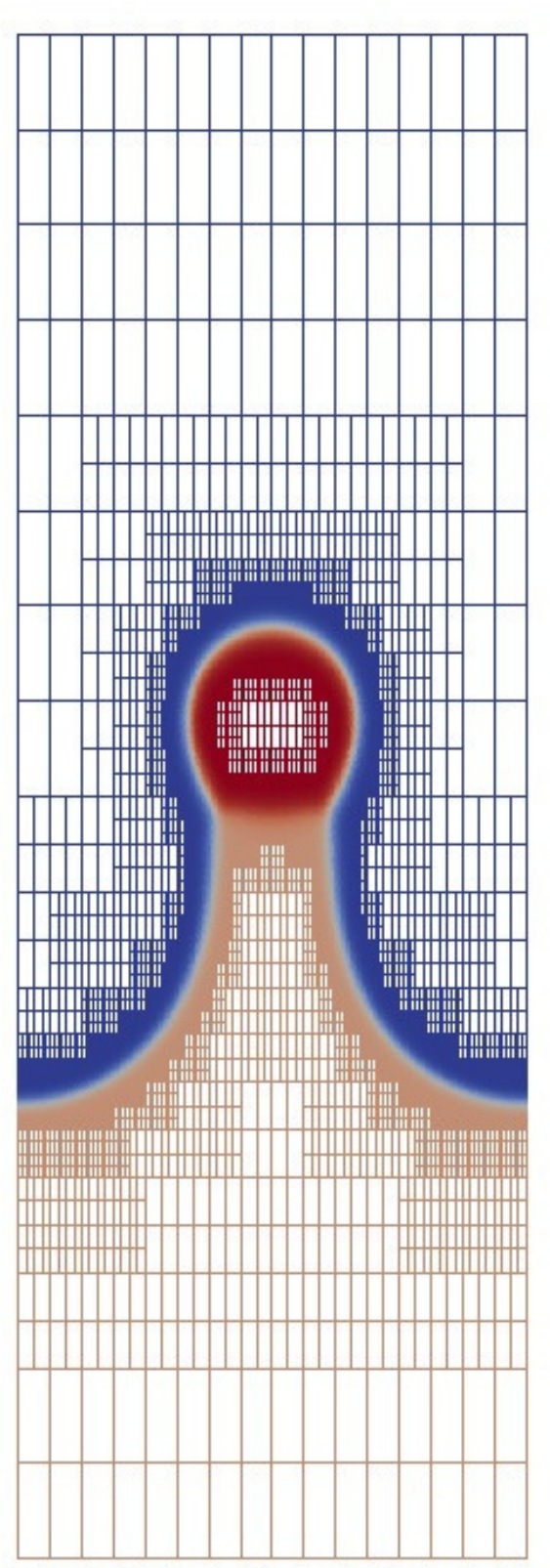}}

  \subfigure[\label{fig:bubble4}]{\includegraphics[width=.17\textwidth]
    {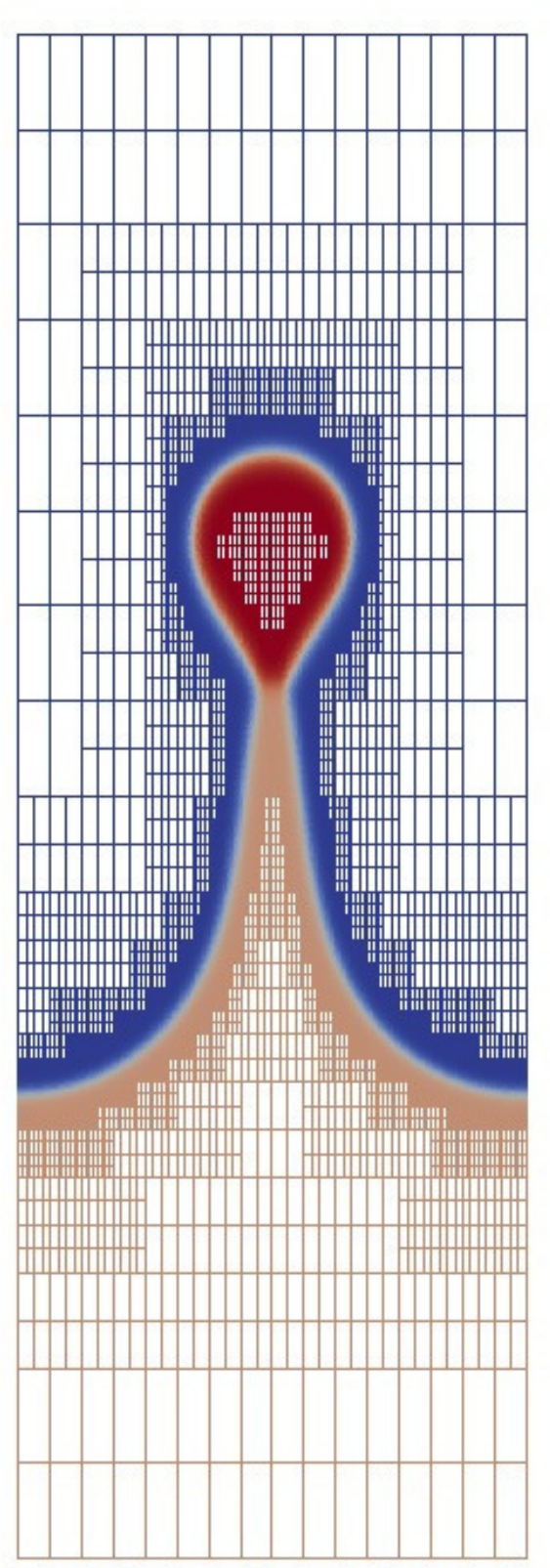}}
  \subfigure[\label{fig:bubble5}]{\includegraphics[width=.17\textwidth]
    {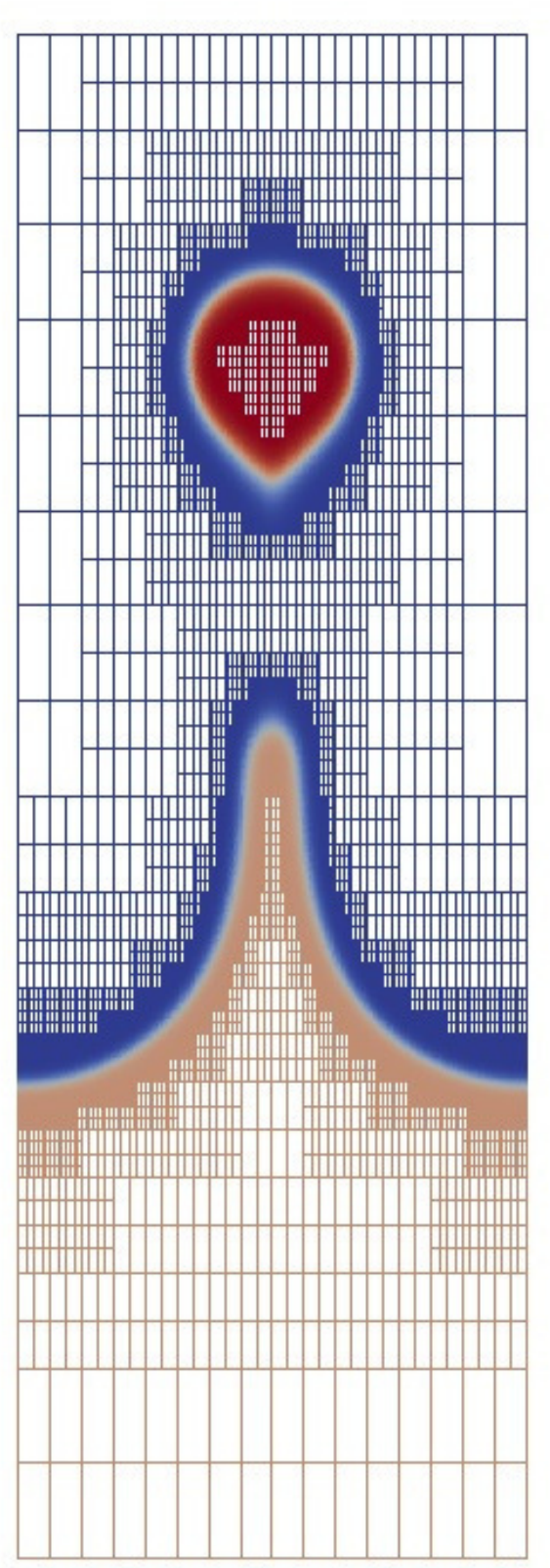}}
  \subfigure[\label{fig:bubble6}]{\includegraphics[width=.17\textwidth]
    {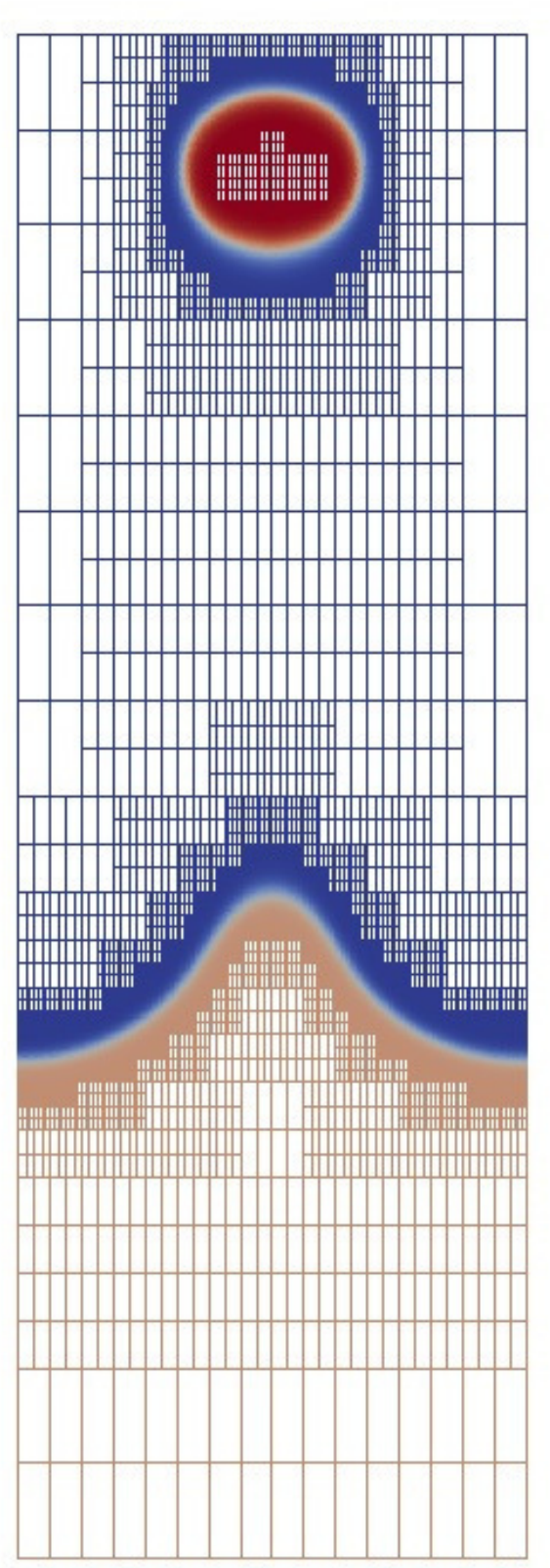}}
  \subfigure[\label{fig:bubble7}]{\includegraphics[width=.17\textwidth]
    {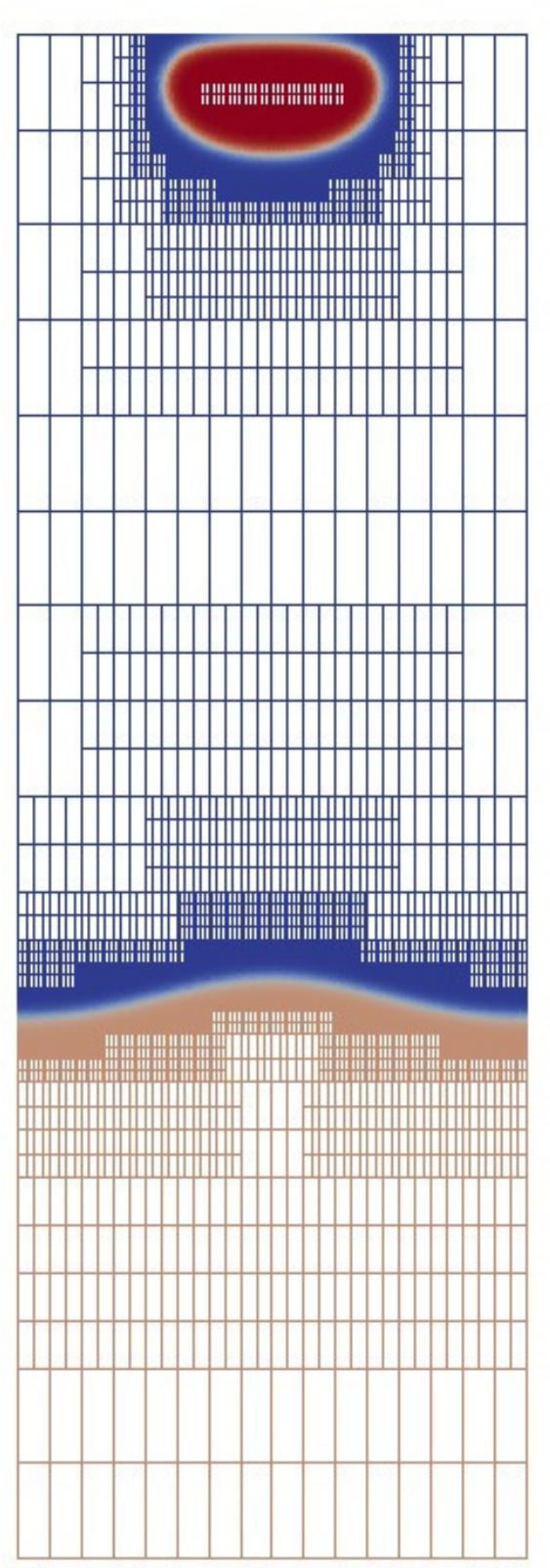}}
\end{center}
  \caption{
    Rise of a fluid bubble ($\psi=1$) that penetrates
    a fluid-fluid interface ($\phi=0$ and $\psi=-1$).
    The computational domain is $[0,1]\times[0,3]$ in all the eight subplots.
    The adaptive mesh refinement and coarsening 
    for the Kelly error estimator are plotted as well.
    \label{fig:bubble}
  }
\end{figure}
We discretize the Stokes system with Taylor-Hood finite elements 
\cite{BrFo91}, a well-known stable pair 
of elements for the velocity and pressure unknowns.
The discretization then becomes 

\begin{eqnarray*}
  && \biggr\langle \,  \left( 1 + \frac{\dt M_1 \gamma_1}{\veps_1^2} 
  (3 (\phiold)^2 - 1) \left(\frac{\psi_h^n-1}{2}\right)^2\right)\phinew 
  \,,\,\eta\,\biggr\rangle \: +  \: \biggr\langle\,\dt \, M_1 \, \gamma_1 \left[
    \left(\frac{\psi_h^n-1}{2}\right)^2 \nabla\phi_{n+1} \right] , \nabla\eta
  \,\biggr\rangle \\
  &&+ \:
  \biggr\langle
  \dt \, \left(\bu_h^n \cdot \nabla \phinew \right), \eta
  \biggr\rangle
  \nonumber 
  = \: \biggr\langle\,\phi_n - 
  \frac{\dt  M_1 \gamma_1}{\veps_1^2} 
  \left( \frac{\psi_h^n-1}{2} \right)^2 (2 (\phiold)^3)
  \eta \,\biggr\rangle, 
\end{eqnarray*}
\begin{eqnarray*}
  &&\biggr\langle \left[ 1 + \frac{\dt M_2\gamma_2}{\veps^2} 
  (3 \left(\psiold\right)^2 - 1) + \frac{1}{2}\dt M_2 \gamma_1
  \left( \frac{1}{2}|\nabla\phi|^2 
  + \frac{1}{4\veps_1^2}(\phi^2-1)^2 \right) \right]
  \psinew  \,,\, \eta \biggr\rangle \\
  &&+ \: \biggr\langle\,
  \dt \, M_2 \gamma_2 \nabla \psi_{n+1}\, ,
  \, \nabla\eta \,\biggr\rangle
  \quad + \quad \biggr\langle \bu_h^n \cdot \psinew, \eta \biggr\rangle
  \nonumber \\
  &&= \:  \biggr\langle 
  \psi_h^n + \frac{\dt M_2 \gamma_2 }{\veps_1^2} (2 \left(\psiold\right)^3) +
  \frac{\dt M_2 \gamma_1}{2} \left[ 
    \frac{1}{2} |\nabla \phi_h^n|^2 + \frac{1}{4 \veps_1^2} 
    \left((\phi_h^n)^2-1\right)^2
    \right], \eta \,\biggr\rangle,
\end{eqnarray*}
\begin{eqnarray*}
  && \mu\biggr\langle\, \nabla \bu_h^n\, , \,\nabla\bw \,\biggr\rangle
  + \biggr\langle \nabla p_h^n, \bw \biggr\rangle
  = \biggr\langle \, \lambda \,\biggr( \Big( \frac{\psi_h^n-1}{2}
  \Big)^2 \, \nabla \phi_h^n \otimes \nabla \phi_h^n + \nabla \psi_h^n
  \otimes \nabla \psi_h^n
  \biggr)\, ,\, \nabla\bw\,\biggr \rangle
  + \biggr \langle \bbf_{ext},\, \bw  \biggr \rangle,
\end{eqnarray*}
and 
\begin{eqnarray*}
  &&\Big\langle \, \nabla \cdot \bu_h^n\,,\, q\,\Big\rangle  = 0.
\end{eqnarray*}
where $\bw$ is a {\bf P2} test function and $q$ and $\eta$ 
are {\bf P1} test functions.

Permuting the equations such that the pressure unknowns appear after the velocities, the
stiffness matrix for Stokes equations has the following block structure
\[
\left(
\begin{matrix}
A & B^T\\
B & 0
\end{matrix}
\right)
 \,
\left(
\begin{matrix}
u\\
p
\end{matrix}
\right)
 =
\left(
\begin{matrix}
f\\
g
\end{matrix}
\right).
\]
A step of block Gaussian Elimination then gives
\begin{eqnarray*}
B A^{-1} B^T p  &=& B A^{-1} f - g \label{psc_1} \: ,  \\
A u &=& f - B^T p \: .
\end{eqnarray*}

The Pressure Schur Complement of the system
$S = B A^{-1} B^T$ ( \cite{Tu99}) then plays the central role
in the linear algebra.  Fortunately, the discrete Laplacian operator 
$A$ is symmetric and positive definite and $B$
has full row rank, an observation that leads to a variety of effective precconditioners 
for this system. 
Our basic solver consists of the following block preconditioner for GMRES iterations:
\[
    P = \left(
    \begin{matrix}
      A & 0\\
      B & -S
    \end{matrix}
    \right), 
  \] or equivalently 
    \[
    P^{-1} = \left(
    \begin{matrix}
      A^{-1} & 0\\
      S^{-1} B A^{-1} & -S^{-1}
    \end{matrix}
    \right) , 
    \]
      {\color{black} so that }
      \[
      P^{-1} \left(
      \begin{matrix}
        A & B^T\\
        B & 0
      \end{matrix}
      \right)
      = 
      \left(
      \begin{matrix}
        I & A^{-1} B^T\\
        0 & I
      \end{matrix}
      \right) .
      \]

We note that this preconditioner is attractive as it reduces the task of inversion of an indefinite system to that 
of solving symmetric and positive definite systems with smaller problem sizes.  Of course, direct inversion 
of $A$ is impractical and thus we use a mass matrix $M_p$ instead in defining the preconditioner. 

The parameters of our tests are set as follows: 
$\dt = 0.01$, $\gamma_1 = \gamma_2 = 0.1$, $\bg =
\left(\,0, 40\,\right)^T$, $\lambda = 0.01$, $M_1 = M_2 = 0.001$, 
$\rho_\bubble -\rho_0=1$, and $\rho_1 = \rho_2 = \rho_0$. 
The gravitational force $\bg$ 
is assigned a large value because the Allan-Cahn
dynamics of our model do not preserve the volumes of the phases. 
Generally, the diffusion
effect will eventually cause the bubble to shrink over time.  This shrinking 
effect 
can be controlled numerically by increasing $\bg$ (\cite{LiSh02}), i.e., 
by reducing the elapsed time in the numerical experiment.

Figure \ref{fig:bubble} shows the upward motion of a fluid bubble ($\psi=1$)
that penetrates a fluid-fluid interface ($\phi=0$ and $\psi=-1$).
The computational domain is $[0,1]\times[0,3]$ in all the eight subplots.
Subplot \ref{fig:bubble0} shows the initial configuration, where the bubble
is immersed in the lower fluid, the fluid-fluid interface is
horizontal, and all the fluids are at rest.
Subplot \ref{fig:bubble1} shows that the fluid-fluid interface
is displayed by the approaching bubble in drift motion.
Accompanying the continuous rise of the bubble, 
the interface breaks to two pieces,
with two three-phase contact lines  formed at the surface of the bubble,
as shown in subplot \ref{fig:bubble2}.
These two three-phase contact  lines then move downward relative to the surface
of the rising bubble, as shown in subplots \ref{fig:bubble3} and 
\ref{fig:bubble4}.
It is observed that in this stage, the bubble exhibits appreciable
deformation, a manifestation of viscoelasticity that arises from
a balance between viscous and capillary forces.
Subplots \ref{fig:bubble4} and \ref{fig:bubble5} 
show that the two three-phase contact  lines merge as they
both reach the stretched bottom of the bubble. Consequently,
the fluid-fluid interfaces, once separated by the intervening bubble,
are joined to form one interface, which is immediately detached from
the bubble. After this pinch off, the bubble continues to rise
in the upper fluid, with a reduced deformation due to
the detachment of the three-phase contact lines from the bubble surface,
as shown in subplots \ref{fig:bubble5} and \ref{fig:bubble6}.
Meanwhile, the surface tension
of the fluid-fluid interface drives it toward the initial (horizontal)
profile. In the last subplot \ref{fig:bubble7}, 
the fluid-fluid interface is already
very close to its initial profile, while the bubble,
with its upward motion stopped by the impermeable boundary,
shows an expected deformation.


\subsection{Investigation of slip}
\begin{figure}[ht!]
\centerline{\includegraphics[width=.4\textwidth]{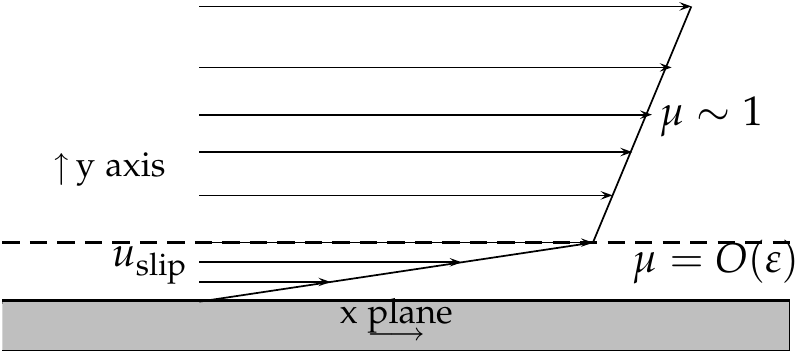}}
  \caption{
    A schematic illustration for fluid slipping modeled by
    a fast variation of tangential velocity across a thin layer
    (diffuse interface) with a small viscosity.
    \label{fig:slip layer}
  }
\end{figure}

As a last application, we incorporate 
slip effects into a simulation. Traditionally, the slip phenomena has been accounted for by
specifying certain boundary conditions. For example, the following
Stokes equations with the so-called Navier boundary condition
\cite{La32}
\begin{eqnarray}
  -\mu \,\Delta \bu + \nabla p &=& 0 \qquad \: \: \:
  \mbox{ in } \Omega \: ,  \label{eq:bvp1}\\
  2 \mu \eps(\bu)_{n\tau} &=& \beta u_\tau^{slip} 
  \quad \mbox{ on } \pa\Omega \: ,  \label{eq:bvp2}
\end{eqnarray}
can be formulated as the following variational problem
\cite{QiWaSh06}
\begin{equation}\label{eq:slip energy}
\min_{\bu, p} \left\{ \int_\Omega 2 \mu \eps(\bu):\eps(\bu) dx
+ \int_{\pa \Omega} \beta \left( u_\tau^{slip} \right)^2 dS \right\} , 
\end{equation}
where $\eps(\bu) = \frac{\nabla \bu + (\nabla \bu)^T}{2}$ is the
symmetric part of the velocity gradient.  Hence, we can see that
the slip boundary conditions contribute to the
dissipation functional.
In our diffuse interface model, we assume slip is induced by  a thin layer (a
diffuse interface associated
with one of the labeling functions, say $\psi$)
of nearly inviscid fluid that surrounds a solid, as depicted in
Figure~\ref{fig:slip layer}.

Formally, we can show that our assumption of a thin
layer on a flat plane can approximate the boundary value problem 
given by (\ref{eq:bvp1}) and
(\ref{eq:bvp2}) in terms of the associated energy as follows.
The bulk energy dissipation functional is given by 
\[
\int_\Omega 2 \mu \eps(\bu):\eps(\bu) dx.
\]
We assume that in the thin layer illustrated in the picture
we have the following velocity profile
\begin{eqnarray*}
  u_{\tau}(x,y) &=& \frac{y}{\veps} u_\tau^{slip}(x),\\
  u_n (x,y) &=& 0,
\end{eqnarray*}
and that the viscosity is proportional to the interfacial
width $\veps$:
\[
\mu = \beta \veps 
\]
in the  thin layer.
Then the energy dissipation in the thin layer reduces to
\begin{eqnarray*}
  \int_{x plane} \int_0^\veps  \beta \veps 
  \left(\frac{u_\tau^{slip}(x)}{\veps}\right)^2
  dy\, dS(x) &=& \int_{x plane} \beta (u_\tau^{slip})^2 dS(x)
\end{eqnarray*}
which is the same as the slip contribution in (\ref{eq:slip energy}).
Therefore, asymptotically (with respect to $\veps$) the diffuse interface model 
with small viscosity inside the interface approximates the traditional sharp interface
model energetically.

We next consider the the effects of our diffuse interface model for slip in an 
application of a solid ball dropping
into a two-phase Stokesian flow.
We distinguish the solid and fluid by $\psi$ and the
fluids by $\phi$. The solid behavior
is achieved by assuming a large bulk viscosity  in the solid ball.
In our implementation, we take
\[
\mu(\psi) = \left\{
\begin{array}{ccc}
30 & \quad & \psi = 1  \\
1 & \quad & \psi = -1   \\
\veps_2 & \quad & |\psi| < 1 
\end{array}
\right. , 
\]
where the larger value of the viscosity ($\mu=30$) is used to model the solid phase
and the smaller ($\mu=\veps_2$) is used to produce the slip effect.


Mathematically, our assumption on viscosity only changes the PDEs of the model 
slightly to the tensor formulation of the Stokes equations:
\[
2 \nabla \cdot ( \mu(\psi)\eps(\bu) ) + \nabla p = \nabla
\cdot \bsigma^e(\phi, \nabla \phi,
\psi, \nabla \psi).
\]
Numerically, we can thus continue to 
use an AMG preconditioner based on the stiffness matrix 
$\langle \mu(\psi) \nabla \bu, \nabla \bw  \rangle$, which is 
spectrally equivalent to the tensor formulation
$\langle 2 \mu(\psi) \eps(\bu), \eps(\bw) \rangle$.

\begin{figure}[ht!]
\begin{center}  \subfigure[\label{fig:slip0}]{\includegraphics[width=.23\textwidth]{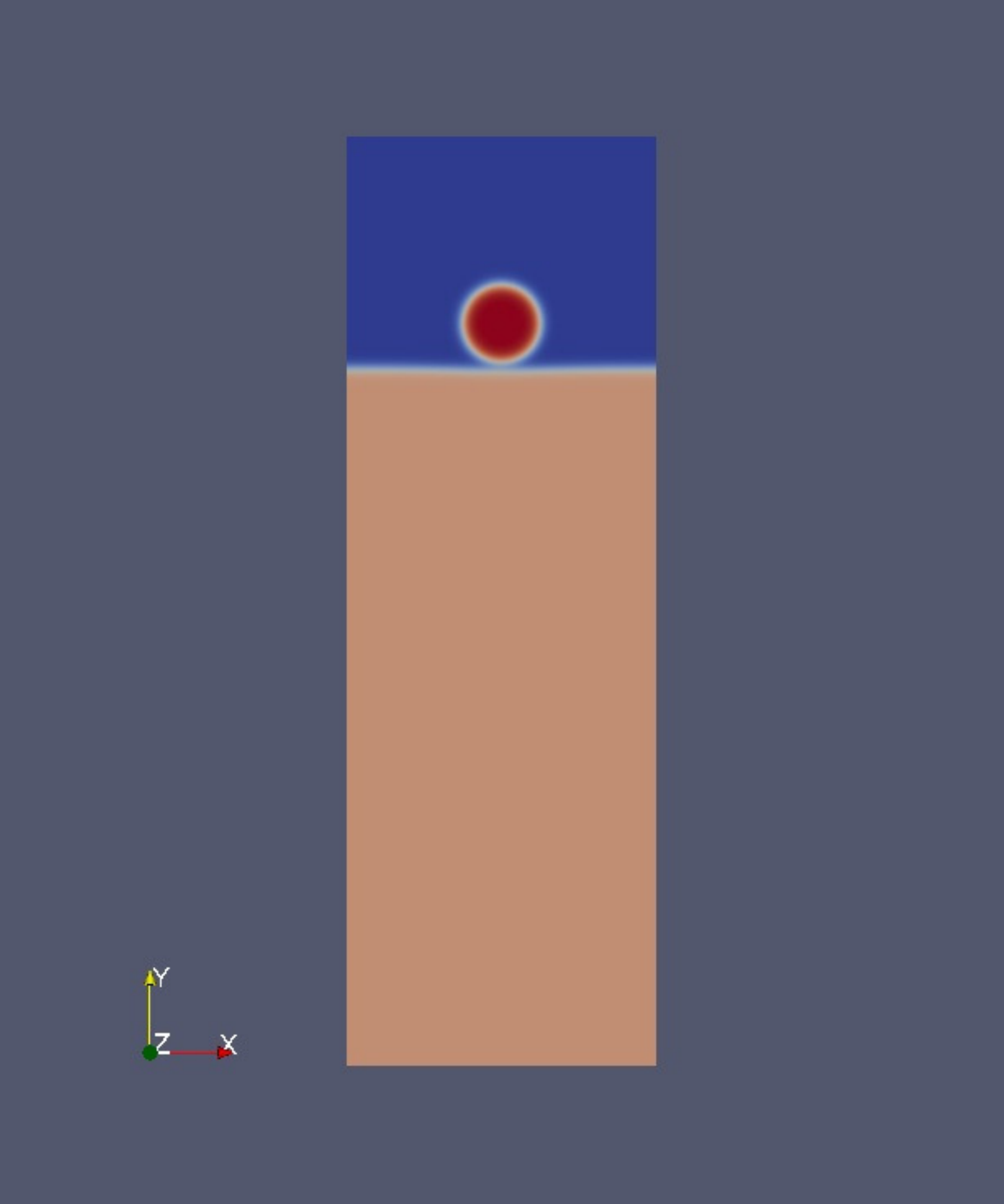}}
  \subfigure[\label{fig:slip1}]{\includegraphics[width=.23\textwidth]{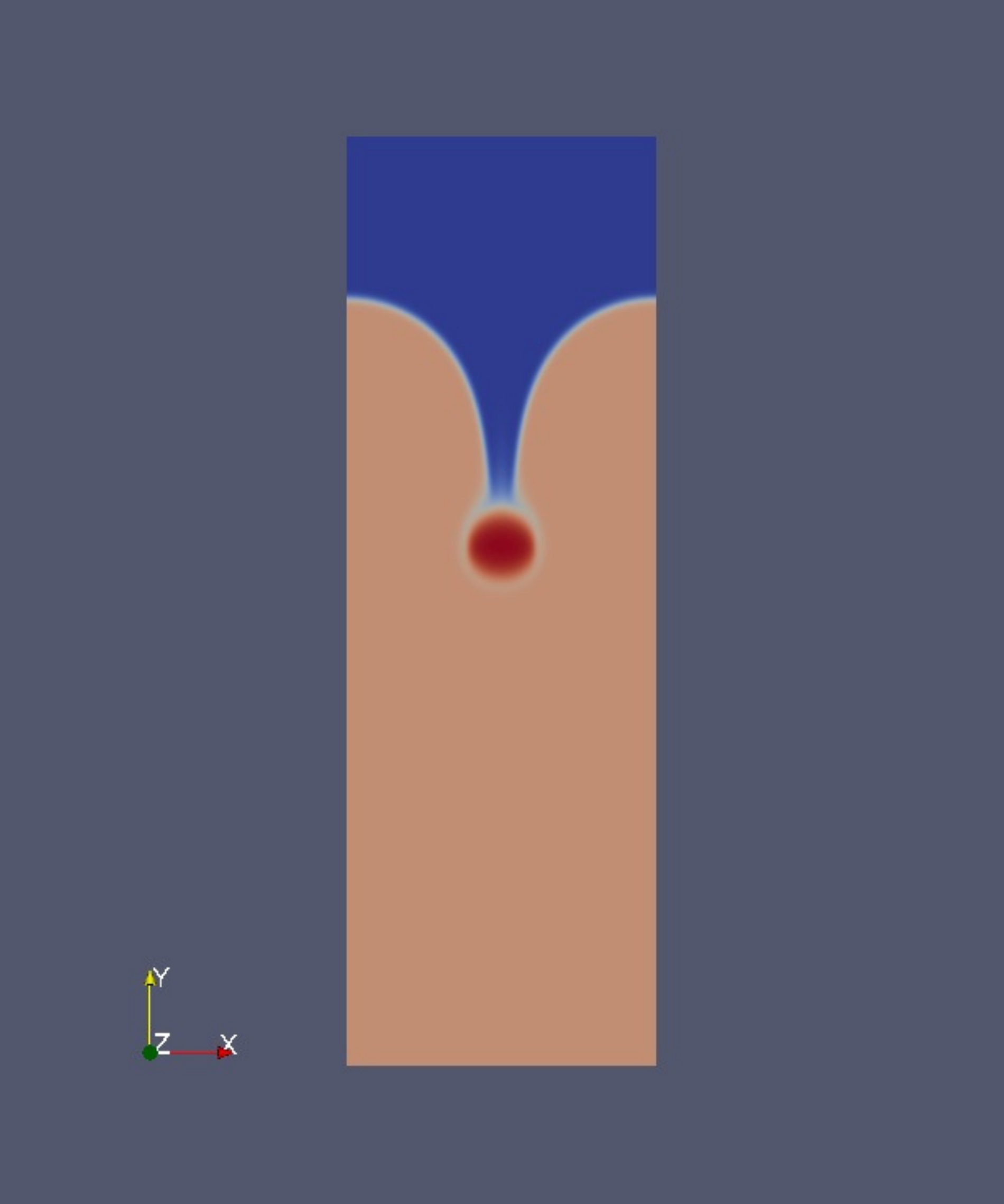}} \hfill \\
  \subfigure[\label{fig:slip2}]{\includegraphics[width=.23\textwidth]{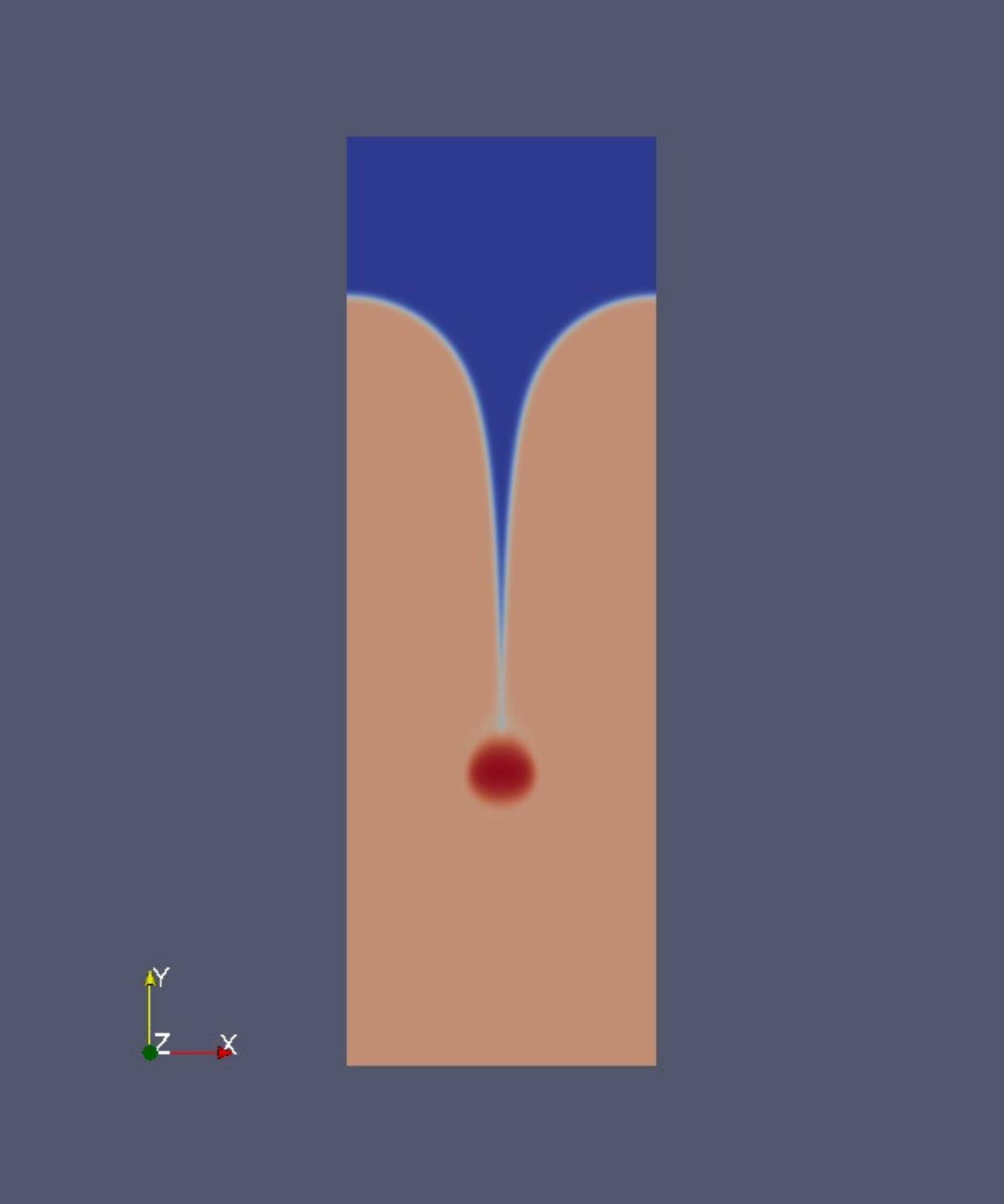}}
  \subfigure[\label{fig:slip3}]{\includegraphics[width=.23\textwidth]{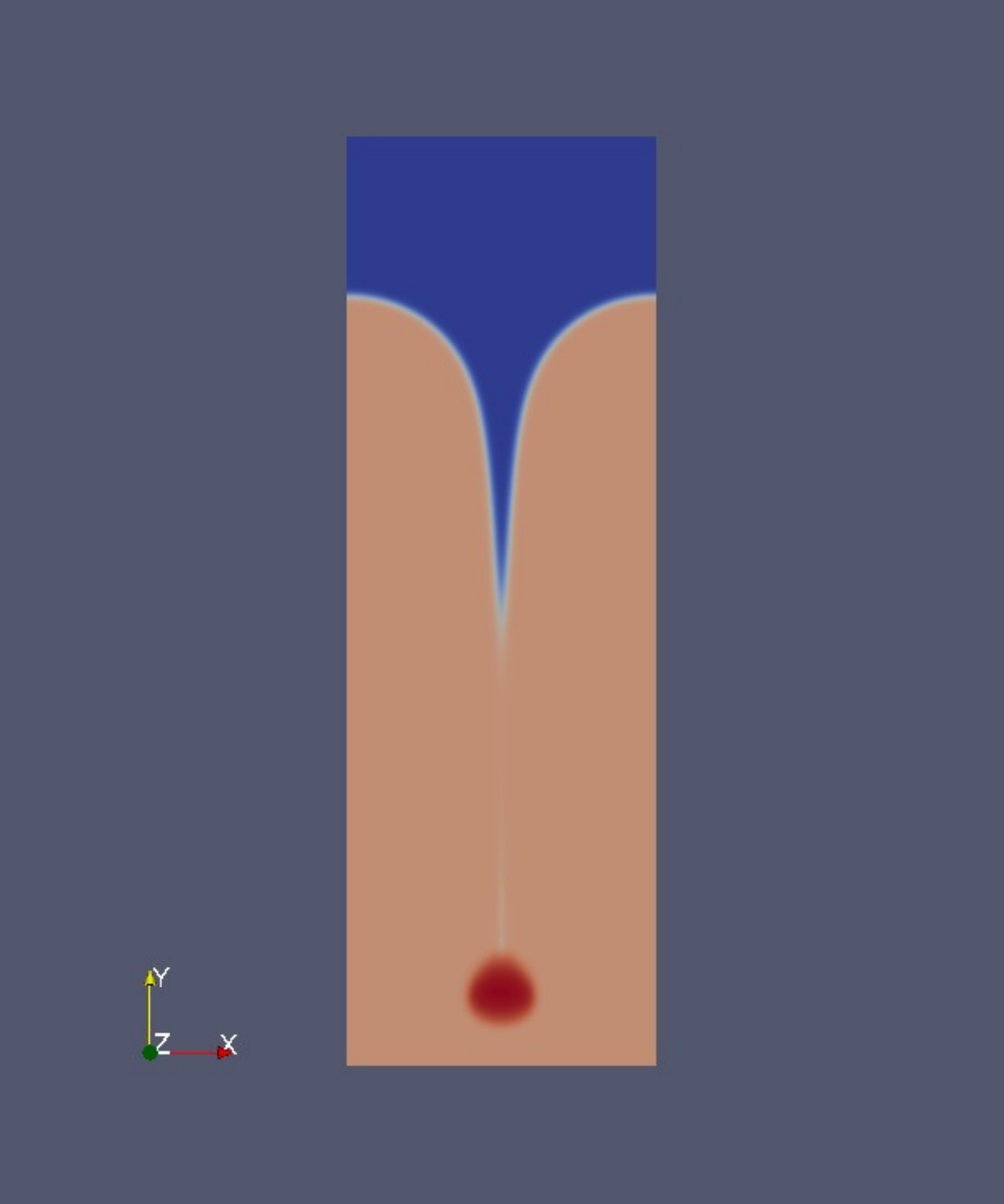}}\hfill
  \end{center}
  \caption{
    A solid particle with a slippery surface falling in a binary fluid.
    Contact lines are formed upon the impact of the particle on
    the fluid-fluid interface. Accompanying the fall of the particle,
    the contact lines quickly move upward relative to the particle surface.
    Upon their arrival at the top of the particle, there is a pinch-off.
    \label{fig:slip}
  }
\end{figure}
\begin{figure}[ht!]
\begin{center}
  \subfigure[\label{fig:noslip0}]{\includegraphics[width=.3\textwidth]{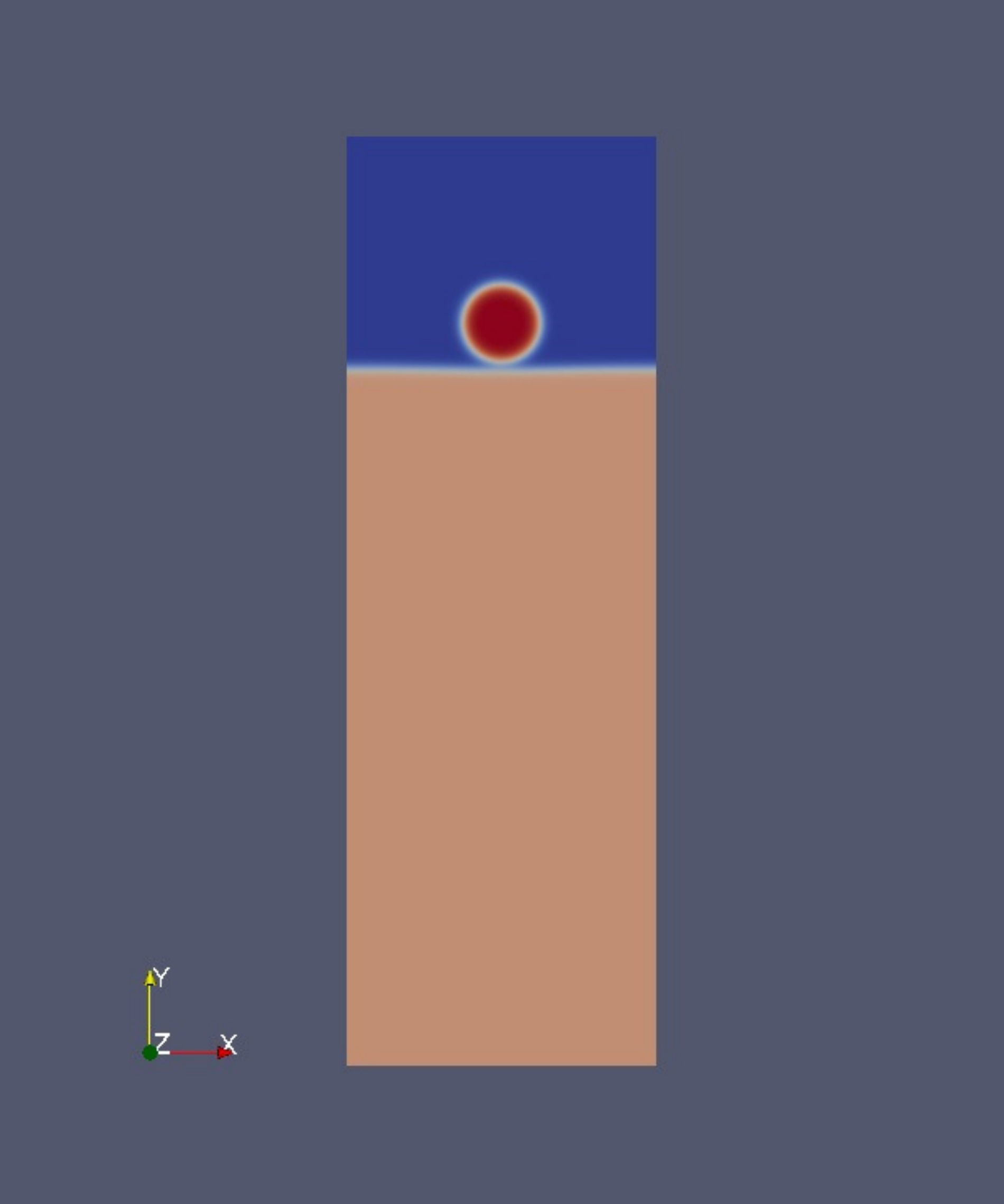}}
  \subfigure[\label{fig:noslip1}]{\includegraphics[width=.3\textwidth]{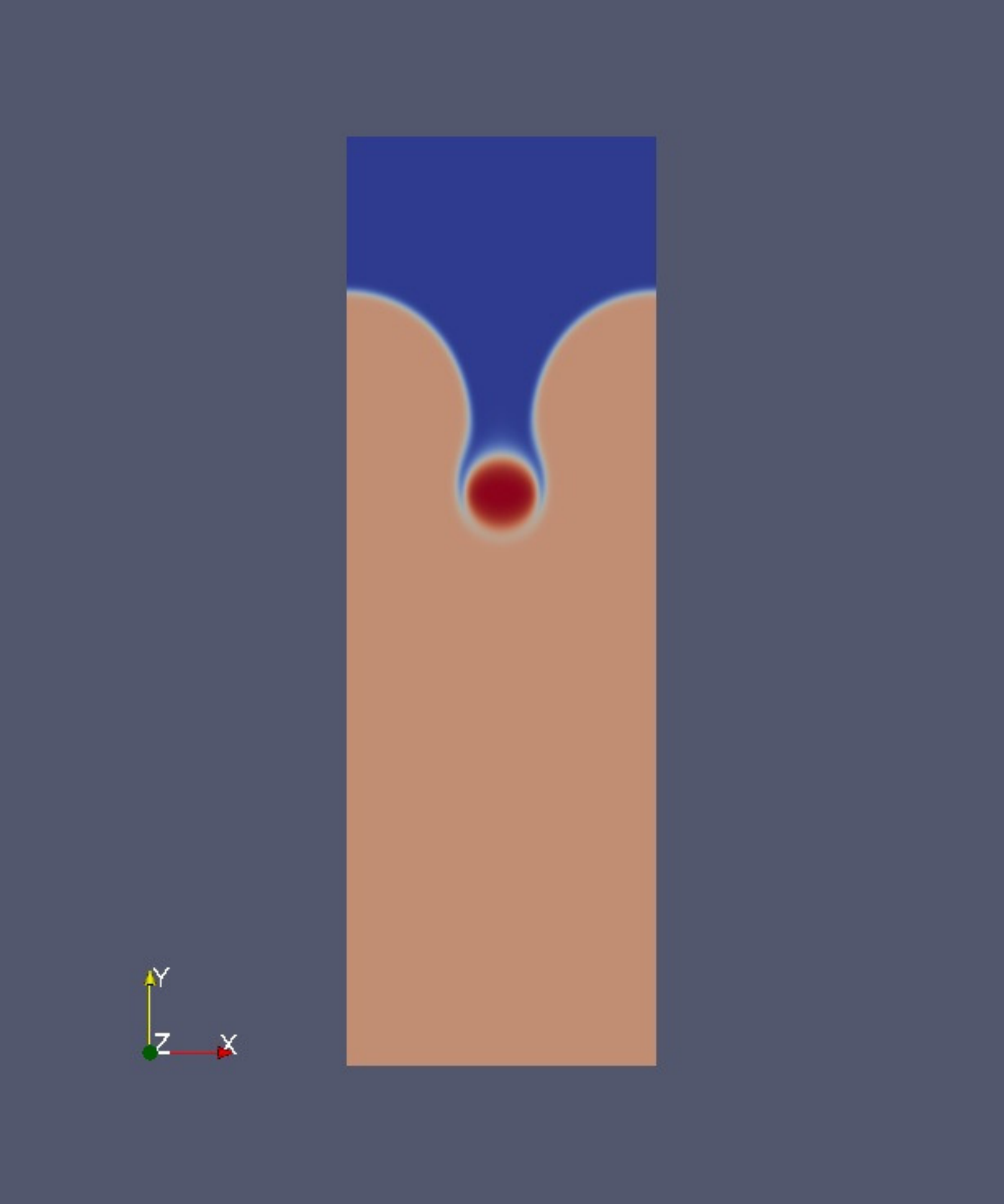}} \\
  \subfigure[\label{fig:noslip2}]{\includegraphics[width=.3\textwidth]{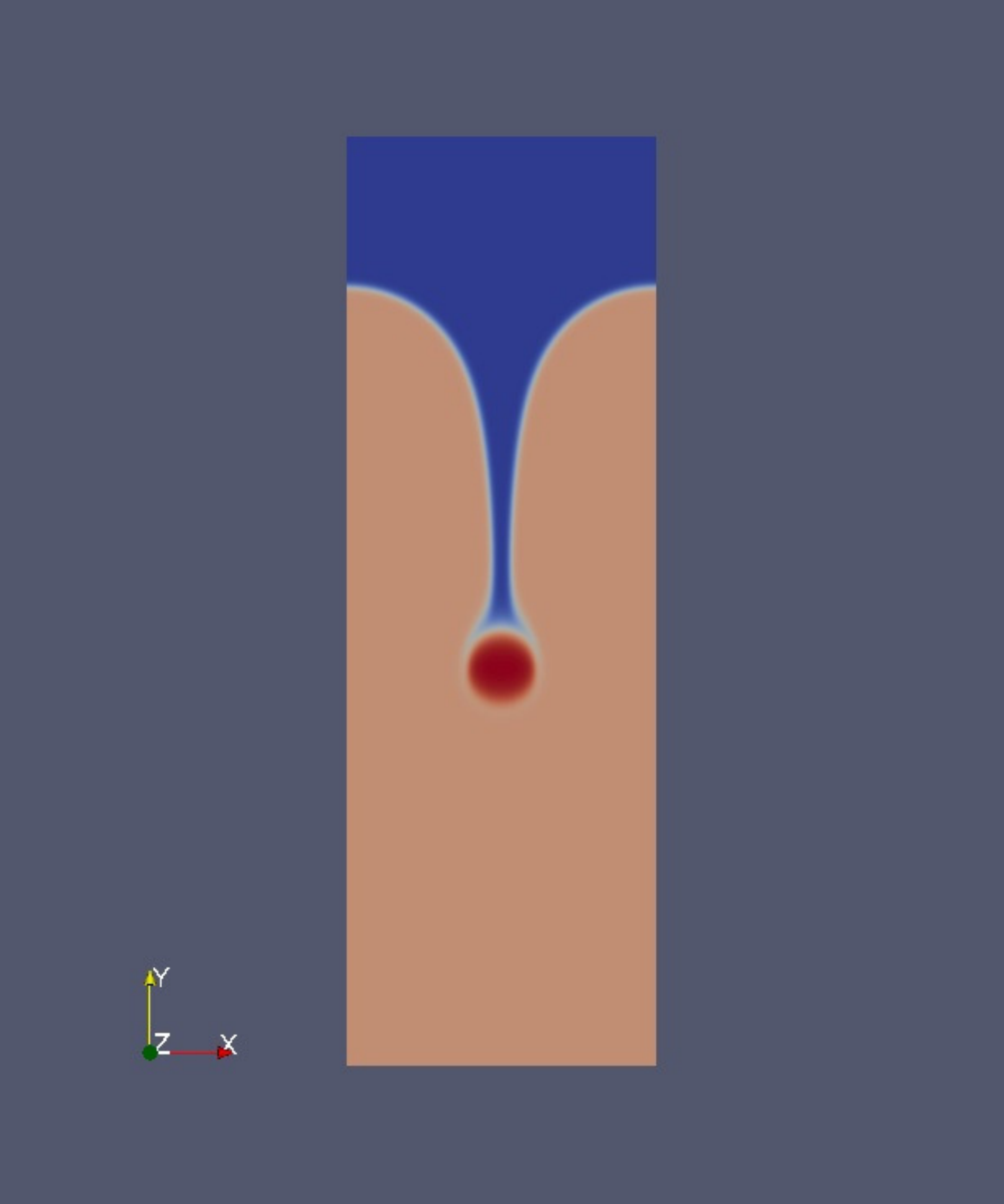}}
  \subfigure[\label{fig:noslip3}]{\includegraphics[width=.3\textwidth]{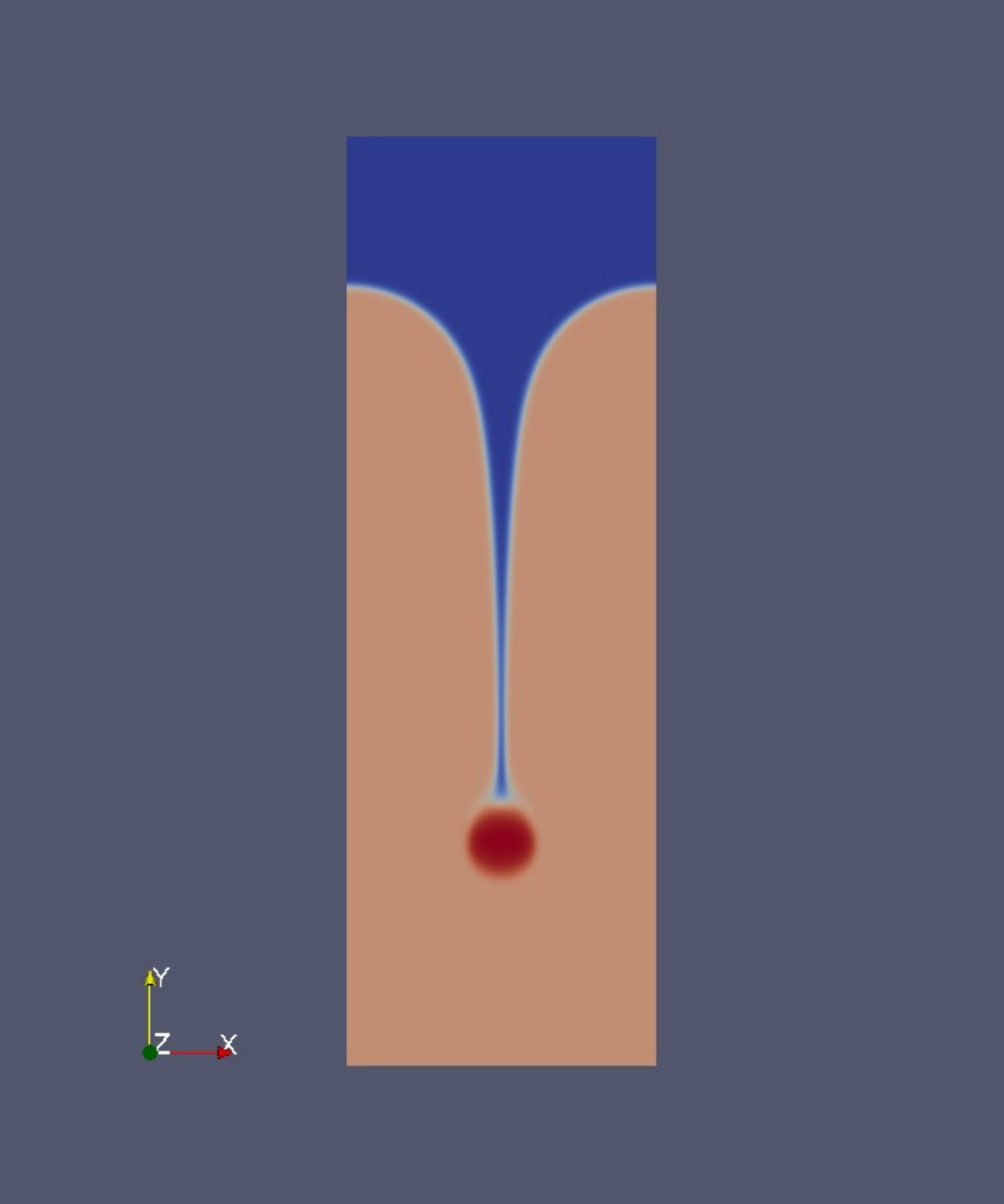}}
    \end{center}
  \caption{
    A solid particle with a non-slippery surface falling in a binary fluid.
    Compared to Figure 5, here the contact line motion is relatively slow.
    Consequently, the particle is wrapped by more upper fluid
    in the early stage of penetration and the pinch-off occurs
    at a larger depth (not shown).
    \label{fig:noslip}
  }
\end{figure}

Numerical simulations have been carried out for a solid particle
that is falling in a two-phase Stokesian flow with and without slip
at the solid surface. The results for a slippery particle are shown
in Figure~\ref{fig:slip} and those for a non-slippery particle are in 
Figure~\ref{fig:noslip}.
It is clearly observed that a slippery surface leads to a more rapid
``pinch-off'' of the upper fluid from the solid particle. Physically,
two contact lines are formed as soon as the falling particle touches
the fluid-fluid interface. Accompanying the fall of the particle,
the contact lines gradually move upward relative to the solid surface.
According to the Onsager principle of least energy dissipation 
\cite{QiWaSh06},
fluid slip would facilitate the contact line motion relative to
the solid surface. The faster moving contact lines at the slippery surface
then lead to an earlier arrival at the top of the particle and
consequently a more rapid pinch-off. This also explains the observation
that during the moment of penetrating the fluid-fluid interface,
there is more upper fluid wrapping the non-slippery particle
with  slower moving contact lines. To summarize, these observations
are consistent with the Onsager principle of least energy dissipation
as fluid slip provides a mechanism to reduce the total rate of
dissipation.

\section{Concluding remarks}\label{sec:Conclude}
We introduced a diffuse interface model to describe
the three-phase dynamics using two phase field variables. The model can be
derived through a variational approach to both the conservative and
the dissipative parts of the dynamics. The applicability of the model
has been demonstrated through simulations for (1) the force balance
at the three-phase contact  line in equilibrium, (2) a rising bubble penetrating
a fluid-fluid interface, and (3) a solid particle falling in
a binary fluid, with fluid slip at solid surface taken into account.
An interesting application of the present model is to further
investigate the effect of particle surface wettability on the impact
of the solid particle on a fluid-fluid interface 
\cite{DYCB07}.
For this purpose, the free energy functional needs to be generalized
to model the different wettabilities of the two fluid phases
on the third (solid) phase.
Work in this direction is currently underway.

The applicability of the model was demonstrated by numerical results that focused on simplified two dimensional models of various multi-phase materials.  Future work in this direction concerns
the development of a three-dimensional parallel code.  This in turn requires further improvements to the numerical models and the numerical algorithms used for approximating the coupled nonlinear systems
of PDEs.


{\flushleft \bf Acknowledgements.} Brannick's work was partially supported by grant NSF-DMS 1320608. Liu's work was partially   
supported by grant NSF DMS-1109107. Sun would like to
thank  the Institute for Mathematics and Its 
Applications (IMA) and Hong Kong University of Science and Technology
as a considerable amount of his work is done during his visits to
the two institutes.

\bibliographystyle{plain}
\bibliography{bibnew}

\def\cprime{$'$}
\begin{thebibliography}{10}

\bibitem{BrFo91}
F.~Brezzi and M.~Fortin.
\newblock {\em Mixed and Hybrid Finite Element Methods}.
\newblock Number~15 in Computational Mathematics. Springer--Verlag, 1991.

\bibitem{BGS96}
Lia Bronsard, Changfeng Gui, and Michelle Schatzman.
\newblock A three-layered minimizer in $r^2$ for a varational problem with a
  symmetric three-well potential.
\newblock {\em Comm. Pure and Applied Math.}, XLIX(673), 1996.

\bibitem{CaHi58}
J.~W. Cahn and J.~E. Hillard.
\newblock Free energy of a nonuniform system. {I}. {I}nterfacial free energy.
\newblock {\em J. Chem. Phys.}, 28:258--267, 1958.

\bibitem{DoEd86}
M.~Doi and S.~F. Edwards.
\newblock {\em The Theory of Polymer Dynamics}.
\newblock Oxford Science Publication, 1986.

\bibitem{DuLiRhWa05a}
Q.~Du, C.~Liu, R.~Ryham, and X.~Wang.
\newblock Modeling the spontaneous curvature effects in static cell membrane
  deformations by a phase field formulation.
\newblock {\em Communications on Pure and Applied Analysis}, 4:537 -- 548,
  2005.

\bibitem{DuLiRyWa05}
Q.~Du, C.~Liu, R.~Ryham, and X.~Wang.
\newblock The phase field formulation of the willmore problem.
\newblock {\em Nonlinearity}, 18:1249 -- 1267, 2005.

\bibitem{DuLiWa04}
Q.~Du, C.~Liu, and X.~Wang.
\newblock A phase field approach in the numerical study of the elastic bending
  energy for vesicle membranes.
\newblock {\em Journal of Computational Physics}, 198(2):450--468, 2004.

\bibitem{DYCB07}
C.~Duez, C.~Ybert, C.~Clanet, and L.~Bocquet.
\newblock Making a splash with water repellency.
\newblock {\em Nature Physics}, 3(180), 2007.

\bibitem{HuReRu94}
W.~Huang, Y.~Ren, and R.~D. Russell.
\newblock Moving mesh partial differential equations (mmpdes) based on the
  equidistribution principle.
\newblock {\em SIAM J. Numer. Anal.}, 31:709--730, 1994.

\bibitem{LoKi05}
J.S.~Lowengrub J.S.~Kim.
\newblock Phase field modeling and simulation of three-phase flows.
\newblock {\em Interfaces Free Bound}, 7:435--466, 2005.

\bibitem{KGZB83}
D.~W. Kelly, J.P. Gago, O.C. Zienkiewicz, and I.Babuska.
\newblock A posteriori error analysis and adaptive proces in the finite element
  method: part i -- error analysis.
\newblock {\em International Journal for Numerical Methods in Engineering},
  19:1593--1619, 1983.

\bibitem{Ki07}
Junseok Kim.
\newblock Phase field computations for ternary fluid flows phase field
  computations for ternary fluid flows.
\newblock {\em Computer Methods in Applied Mechanics and Engineering},
  196(45):4779--4788, 2007.

\bibitem{La32}
H.~Lamb.
\newblock {\em Hydrodynamics}.
\newblock Cambridge, 6th edition, 1932.

\bibitem{LeLiZh08}
Zhen Lei, Chun Liu, and Yi~Zhou.
\newblock Global solutions for incompressible viscoelastic fluids.
\newblock {\em Arch. Rational Mech. Anal.}, 188, 2008.

\bibitem{LiLiZh05}
F.~H. Lin, C.~Liu, and P.~Zhang.
\newblock On a micro-macro model for polymeric fluids near equilibrium.
\newblock {\em Comm. Pure Appl. Math.}, LIX:1--29, 2005.

\bibitem{LiSh02}
C.~Liu and J.~Shen.
\newblock A phase field model for the mixture of two incompressible fluids and
  its approximation by a fourier-spectral method.
\newblock {\em Physica D}, 179:211--228, 2003.

\bibitem{QWS03}
T.~Qian, X.~P. Wang, and P.Sheng.
\newblock Molecular scale contact line hydrodynamics of immiscible flows.
\newblock {\em Phys. Rev. E}, 68(016306), 2003.

\bibitem{QiWaSh06}
T.~Qian, X.~P. Wang, and P.~Sheng.
\newblock A variational approach to the moving contact line hydrodynamics.
\newblock {\em J. Fluid Mech.}, 564:333 -- 360, 2006.

\bibitem{EES83}
S.~S.~Eisenstat, H.~Elman, and M.~Schultz.
\newblock Variational iterative methods for nonsymmetric systems of linear
  equations.
\newblock {\em SIAM J. of Num. Anal.}, 20(2):345--357, 1983.

\bibitem{SY10}
Jie Shen and Xiaofeng Yang.
\newblock A phase-field model and its numerical approximation for two-phase
  incompressible flows with different densities and viscosities.
\newblock {\em SIAM J. Sci. Computing}, 32(1159), 2010.

\bibitem{SL08}
H.~Sun and C.~Liu.
\newblock On energetic variational approaches in modeling the nematic liquid
  crystal flows.
\newblock {\em DCDS-A}, 23:455--475, 2008.

\bibitem{Tu99}
Stefan. Turek.
\newblock {\em Efficient solvers for incompressible flow problems : an
  algorithmic and computational approach}.
\newblock Springer, 1999.

\bibitem{YuFeLiSh05b}
P.~Yue, J.~Feng, C.~Liu, and J.~Shen.
\newblock A diffuse-interface method for simulating two-phase flows of complex
  fluids.
\newblock {\em Journal of Fluid Mechanics}, 515:293--317, 2005.

\end{thebibliography}
\end{document}